\documentclass[twoside,11pt]{article}

\usepackage{blindtext}

\usepackage{jmlr2e}

\usepackage{graphicx}
\usepackage{booktabs} 

\usepackage[flushleft]{threeparttable}

\usepackage{epsf}
\usepackage{fancyhdr}
\usepackage{graphics}
\usepackage{graphicx}
\usepackage{psfrag}
\usepackage{microtype}
\usepackage{verbatim}

\usepackage{color}
\usepackage{amsfonts}
\usepackage{amsmath}
\usepackage{amssymb,bbm}

\usepackage{algorithm}
\usepackage{algorithmic}

\usepackage{caption,subcaption}
\usepackage{sidecap}
\sidecaptionvpos{figure}{c}

\usepackage{url}
\usepackage{hyperref}

\newcommand{\KL}{\mathbf{KL}}
\newcommand{\OT}{\mathbf{OT}}
\newcommand{\UOT}{\mathbf{UOT}}

\long\def\comment#1{}
\comment{
\setlength{\topmargin}{0 in}
\setlength{\textwidth}{5.5 in}
\setlength{\textheight}{8 in}
\setlength{\oddsidemargin}{0.5 in}
}

\newcommand{\norm}[1]{\left\lVert#1\right\rVert}

\newcommand{\Br}{\mathbb{R}}
\newcommand{\one}{\textbf{1}}

\def\Xx{\mathcal{X}}

\def\Yy{\mathcal{Y}}

\def\aA{\mathbf{a}}
\def\bB{\mathbf{b}}
\def\xX{\mathbf{x}}
\def\yY{\mathbf{y}}
\def\uU{\mathbf{u}}
\def\vV{\mathbf{v}}

\def\tT{\mathbf{t}}
\def\gG{\mathbf{g}}

\graphicspath{{./image/}}

\newcommand{\ba}{\begin{array}}
\newcommand{\ea}{\end{array}}

\newcommand{\argminE}{\mathop{\mathrm{argmin}}}   

\newcommand{\at}[2][]{#1|_{#2}}

\newcommand{\alglinelabel}{%
  \addtocounter{ALC@line}{-1}
  \refstepcounter{ALC@line}
  \label
}

\hypersetup{ hidelinks }

\ShortHeadings{Efficient Algorithms and Approximation Error for Unbalanced Optimal Transport}{Q. M. Nguyen, H. H. Nguyen, Y. Zhou and L. M. Nguyen}
\firstpageno{1}

\usepackage{lastpage}
\jmlrheading{24}{2023}{1-\pageref{LastPage}}{10/22; Revised
9/23}{12/23}{22-1158}{Quang Minh Nguyen, Hoang H. Nguyen, Yi Zhou, and Lam M. Nguyen}
\ShortHeadings{Efficient Algorithms and Approximation Error for Unbalanced Optimal Transport}{Q. M. Nguyen, H. H. Nguyen, Y. Zhou and L. M. Nguyen}

\begin{document}

\title{On Unbalanced Optimal Transport: Gradient Methods, Sparsity and Approximation Error}

\author{\name Quang Minh Nguyen \email nmquang@mit.edu \\
       \addr Electrical Engineering and Computer Science\\
  Massachusetts Institute of Technology\\
 Cambridge, MA, USA 
       \AND
       \name Hoang H. Nguyen \email hnguyen455@gatech.edu \\
       \addr Industrial and Systems Engineering\\
   Georgia Institute of Technology\\
    Atlanta, GA, USA
       \AND
       \name Yi Zhou \email yi.zhou@ibm.com  \\
       \addr IBM Research\\
  Almaden Research Center\\
  San Jose, CA, USA
       \AND
       \name Lam M. Nguyen \email LamNguyen.MLTD@ibm.com \\
       \addr IBM Research \\
   Thomas J. Watson Research Center\\
   Yorktown Heights, NY, USA
       }

\editor{Silvia Villa}

\maketitle

\begin{abstract}
We study the Unbalanced Optimal Transport (UOT)  between two measures of possibly different masses with at most $n$ components, where the marginal constraints of  standard Optimal Transport (OT) are relaxed via Kullback-Leibler divergence with regularization factor $\tau$. Although only Sinkhorn-based UOT solvers have been analyzed in the literature with the iteration complexity of ${O}\big(\tfrac{\tau  \log(n)}{\varepsilon}  \log\big(\tfrac{\log(n)}{{\varepsilon}}\big)\big)$ and per-iteration cost of $O(n^2)$  for achieving the desired error $\varepsilon$, their positively dense output transportation plans strongly hinder the practicality. On the other hand, while being vastly used as heuristics for computing UOT in modern deep learning applications and having shown success in sparse OT problem, gradient methods applied to UOT have not been formally studied. In this paper, we propose a novel algorithm based on Gradient Extrapolation Method (GEM-UOT) to find  an $\varepsilon$-approximate solution to the  UOT problem in $O\big(  \kappa  \log\big(\frac{\tau  n}{\varepsilon}\big) \big)$ iterations with $\widetilde{O}(n^2)$ per-iteration cost, where $\kappa$ is the condition number depending on only the two input measures.  Our  proof technique is based on a   novel dual formulation of the squared $\ell_2$-norm  UOT objective, which  fills  the lack of sparse UOT literature  and also leads to a new characterization of approximation error between UOT and OT. To this end, we further present a novel approach of OT retrieval from UOT, which is based on GEM-UOT with fine tuned $\tau$ and a post-process projection step.
Extensive experiments on synthetic and real datasets validate our theories and demonstrate the favorable performance of our methods in practice. We  showcase GEM-UOT on the task of color transfer in terms of both the quality of the transfer image and the sparsity of the transportation plan. 

\end{abstract}


\section{Introduction}\label{sec_intro}
The optimal transport (OT) problem  originated from the need to find the optimal cost to transport masses from one distribution to another distribution \citep{Villani-09}. 
While initially developed by theorists, OT has found widespread applications in statistics and machine learning (ML) (see e.g. \citep{ho2017multilevel,arjovsky2017wasserstein,Rabin2015ConvexCI}). However, standard OT requires the restricted assumption that the input measures are  normalized to unit mass, which facilitates the development of the Unbalanced Optimal Transport (UOT) problem between two measures of possibly different masses \citep{chizat2018scaling}. The class of UOT problem relaxes OT’s marginal constraints. Specifically, it is a regularized version of Kantorovich
formulation placing penalty functions on the marginal
distributions based on some given divergences \citep{Liero_2017}. While there have been several divergences considered by the literature, such as squared $\ell_2$ norm \citep{blondel2018smooth}, $\ell_1$ norm \citep{l1_uot}, or general $\ell_p$ norm \citep{lee2019parallel}, UOT with Kullback-Leiber  (KL) divergence \citep{chizat2018scaling} is the most prominent and has been used in statistics and machine learning \citep{wassloss}, deep learning \citep{Yang_Scalable_2019}, domain adaptation \citep{balaji2020robust, Fatras2021UnbalancedMO}, bioinformatics \citep{Schiebinger_Optimal_2019}, and OT robustness \citep{balaji2020robust, le2021robust}. Throughout this paper, we refer to UOT penalized by KL divergence as simply UOT, unless otherwise specified. We hereby define some notations and formally present our problem of interest. 

\textbf{Notations:}
We let $[n]$ stand for the set $\{1, 2, . . . , n\}$ while $\Br^n_+$ stands for the set of all vectors in $\Br^n$  with nonnegative entries. For a vector $\xX \in R^n$ and $p \in [1, \infty)$, we denote $\|\xX\|_p$ as its $\ell_p$-norm and $diag(\xX) \in \Br^{n\times n}$ as the diagonal matrix with $diag(\xX)_{ii}= x_i$.
Let A and B be two matrices of size $n\times n$, we denote their Frobenius inner product as: $\left\langle A , B\right\rangle = \sum_{i,j=1}^n A_{ij} B_{ij}$. $\one_n$ stands for a vector of length $n$ with all of its components equal to $1$. The KL divergence between two vectors $\xX, \yY \in \Br^n$ is defined as $\KL(x || y) = \sum_{i = 1}^n x_i \log \left(\frac{x_i}{y_i} \right) - x_i + y_i$. The entropy of a matrix $X\in \Br^{n\times n}_+$ is given by $H(X) = -\sum_{i,j=1}^n X_{ij} (\log(X_{ij})-1) $. 

\textbf{Unbalanced Optimal Transport:}
For a couple of finite measures with possibly different
total mass $\aA = (a_1,...,a_n) \in \mathbb{R}^{n}_+$ and  $\bB = (b_1,...,b_n) \in \mathbb{R}^{n}_+$, we denote $\alpha = \sum_{i=1}^n a_i$ and $\beta = \sum_{i=1}^n b_i$ as the total masses, and $a_{min} = \min_{1 \leq i \leq n}\{a_i\}$ and $b_{min} = \min_{1 \leq i \leq n}\{b_i\}$ as the minimum masses. The UOT problem can be written as:
\begin{align}
\label{UOT_func}
    \UOT_{\KL}(\aA, &\bB) =  \min_{X \in \Br_{+}^{n \times n}} \big\{ f(X) := \left\langle C, X\right\rangle  + \tau \KL(X \one_{n} || \aA) + \tau \KL(X^{\top} \one_{n} || \bB) \big\},
\end{align}
where $C$ is a given cost matrix, $X$ is a transportation plan, $\tau > 0$ is a given regularization parameter. When $\aA^{\top} \one_n = \bB^{\top} \one_n$ and $\tau \rightarrow \infty$, \eqref{UOT_func} reduces to a standard OT problem. Let $X_f \in \argminE_{X \in \Br_{+}^{n \times n}} f(X)$ be  optimal transportation plan of the UOT problem \eqref{UOT_func}.


\begin{definition}
For  $\varepsilon > 0$, $X$ is an $\varepsilon$-approximate transportation plan of $\UOT_{\KL}(\aA, \bB)$ if:
\begin{align}
    &f(X)=\langle C, X\rangle + \tau \KL(X \one_{n} || \aA)  
     + \tau \KL(X^{\top} \one_{n} || \bB) \leq \UOT_{\KL}(\aA, \bB) + \varepsilon.
\end{align}


\end{definition}


\subsection{Open Problems and Motivation}
\label{sec:open_problems}
We hereby discuss the open problems within the UOT problem to better facilitate the motivation underpinning our work.

\textbf{Nascent literature of gradient methods for UOT:} While the existing work has well studied the class of OT problems using gradient methods or variations of the Sinkhorn algorithm \citep{guminov2021accelerated, altschuler2017near}, the complexity theory and algorithms for UOT remain nascent despite its recent emergence. \cite{pham2020unbalanced} shows that  Sinkhorn algorithm can solve the UOT problem in ${O}\big(\frac{\tau n^2 \log( n)}{\varepsilon} \log \big(\frac{\log(n)}{\varepsilon}\big)\big)$ up to the desired error $\varepsilon$. \cite{Chapel2021UnbalancedOT} proposes Majorization-Minimization (MM) algorithm, which is  Sinkhorn-based, as a numerical solver with $O(n^{3.27})$ empirical complexity. Thanks to their compatibility, gradient methods have been vastly used as heuristics for computing UOT in  deep learning applications \citep{balaji2020robust, Yang_Scalable_2019} and multi-label learning \citep{wassloss}. Nevertheless, no work has formally studied gradient methods applied to UOT problems. 
Principled study of this topic could lead to preliminary theoretical justification and guide the more refined design of gradient methods for UOT.

\textbf{Lack of sparse UOT literature:} There is currently no study of sparse UOT in the literature, which nullifies its usage in many applications \citep{PITIE2007123, courty2016optimal, muzellec2016tsallis} in which only sparse transportation plan  is of interest. On a side note, entropic regularization imposed by the Sinkhorn algorithm keeps the transportation plan dense and strictly positive \citep{sparseOT_solomon, blondel2018smooth}.

\textbf{Approximation of OT:}
The approximability of standard OT via UOT has remained an open problem. 
Recently, \citep{Chapel2021UnbalancedOT} discusses the possibility to approximate OT transportation plan using that of UOT and empirically verifies it. 
Despite the well known fact that OT is recovered from UOT when $\tau \to \infty$ and the masses are balanced, no work has analyzed the rate under which UOT converges to OT\footnote{For squared $\ell_2$ norm penalized UOT variant, \citep{blondel2018smooth} showed that  error in terms of transport distance is $O(n^2/ \tau)$.}. Characterization of approximation error between UOT and OT could give rise to the   new perspective of OT retrieval from UOT, and facilitate further study of the  UOT geometry and statistical bounds from the lens of OT \citep{otgeo1} in the regime of large $\tau$. 

\textbf{Bottleneck of robust OT computation via UOT:}  The UOT formulation with small $\tau$ has been vastly used as a robust variant of OT\footnote{While \cite{Fatras2021UnbalancedMO} uses exactly UOT as robust OT, \cite{le2021robust} considers UOT with an additional normalization constraint, i.e. $\|X\|_1 = 1$,  which can be solved with usual UOT algorithm with  an extra normalization step.}  \citep{Fatras2021UnbalancedMO, balaji2020robust, le2021robust}.  In this paper, we further highlight a major bottleneck of this approach having been neglected by the literature: though relaxing OT as  UOT admits robustness to outliers, the computed UOT distance far deviates from the original OT distance. Under the naive choice of $\tau = 1$  \citep{Fatras2021UnbalancedMO, balaji2020robust, le2021robust}, we  observe that on CIFAR-10 as an example real dataset the UOT distance differs from OT by order of dozens (see Figure \ref{fig_004}). Such large deviation from the original metric   can be  detrimental to target applications where the  distance itself is of interest  \citep{genevay2018learning}. Another inherent limitation is that the  solution returned by UOT is only a ``relaxed" transportation plan that does not respect the marginal constraints of standard OT. This naturally raises the question whether despite the relaxation of OT as UOT, there is  a way to retrieve the original OT distance and transportation plan.

\subsection{Contributions}


In this paper, we provide a comprehensive study of UOT that addresses all the aforementioned open problems and challenges of UOT literature. We consider the setting where KL divergences are used to penalize the marginal constraints.
Our contributions can be summarized as follows:

\begin{itemize}
    \item  We  provide a novel dual formulation of the squared $\ell_2$-norm regularized UOT objective, which is the basis of our algorithms. This could facilitate further algorithmic development for solving sparse or standard UOT, since the current literature is limited to the dual formulation of entropic-regularized UOT problem \citep{chizat2018scaling}. 
    \item  Based on the Gradient Extrapolation Method (GEM), we propose GEM-UOT algorithm for finding   an $\varepsilon$-approximate solution to the  UOT problem in $O\big(  \kappa  \log\big(\frac{\tau  n}{\varepsilon}\big) \big)$ iterations with $\widetilde{O}(n^2)$ per-iteration cost, where $\kappa$ is the condition number that depends on only the two input measures. 
    Through GEM-UOT, we present the first i) principled study of gradient methods for UOT, ii) sparse UOT solver, and iii)   algorithm that lifts the linear dependence on  $\tau$ of Sinkhorn  \citep{pham2020unbalanced}, a bottleneck in the regime of large $\tau$ \citep{sejourne2022faster}.
    
    \item To the best of our knowledge, we establish the first characterization of the approximation error between UOT and OT (in the context of KL-penalized UOT). In particular, we show that both of UOT's transportation plan and transport distance converge to OT's marginal constraints and transport distance with the rate $O(\frac{poly(n)}{\tau})$. This result  opens up directions that use UOT to approximate standard OT and accentuate the importance of our proposed GEM-UOT, which is the first to achieve logarithmic dependence on $\tau$  in the literature.

    \item Inspired by our results on approximation error, we bring up a new perspective of OT retrieval from UOT with fine-tuned $\tau$. In particular, we present GEM-OT, which obtains an $\varepsilon$-approximate solution to the standard OT problem by performing a post-process projection of UOT solution. GEM-OT is of order $\widetilde{O}(\kappa n^2)$.
    While lots of work have naively relaxed OT into UOT to enjoy its unconstrained optimization structure and robustness to outliers, our results on $\tau$ make a further significant step by equipping such relaxation with provable guarantees.
\end{itemize}

\noindent
\textbf{Paper Organization:} The rest of the paper is organized as follows. We introduce the background of regularized UOT problems and present our novel dual formulation in Section~\ref{sec_background}. In Section~\ref{sec_analysis}, we analyze the complexities our proposed algorithms GEM-UOT. The results on approximability of OT via UOT 
are established in Section ~\ref{sec:approximation}.  In Section~\ref{sec_experiment}, we experiment on both synthetic and real datasets to compare our algorithm with the state-of-the-art Sinkhorn algorithm in terms of both the performance and the induced sparsity, empirically verify our theories on approximation error, and showcase our framework on the task of color transfer. We conclude the paper in Section~\ref{sec_conclusion}.

\section{Background}


\label{sec_background}
 In this Section, we first present the entropic regularized UOT problem, used by Sinkhorn algorithm. Then we consider the squared  $\ell_2$-norm regularized UOT and derive a new dual formulation, which is the basis for our algorithmic design.


\subsection{Entropic Regularized UOT}
Inspired by the literature of the entropic regularized OT problem, the entropic version of the UOT problem has been considered. 
The problem is formulated as:
\begin{equation}
\label{UOT_regu_func}
     \min_{X \in \Br_{+}^{n \times n}} \left\langle C, X\right\rangle - \eta H(X) + \tau \KL(X \one_{n} || \aA) + \tau \KL(X^{\top} \one_{n} || \bB),
\end{equation}
where $\eta > 0$ is a given regularization parameter.
By \citep{chizat2018scaling}, optimizing the Fenchel-Legendre dual of the above entropic regularized UOT is equivalent to:
\begin{equation}
\label{dual_UOT}
    \min_{\uU, \vV \in \Br^{n}}  \eta \sum_{i,j=1}^{n} exp\bigg(\frac{u_i+v_j - C_{ij}}{\eta}\bigg) + \tau \left\langle  e^{-\uU/\tau} , \aA\right\rangle + \tau \left\langle  e^{-\vV/\tau} , \bB\right\rangle.
\end{equation}

\begin{remark}
\label{entropic_downvote}
We  note that existing algorithm   \citep{pham2020unbalanced}  must set small $\eta = O(\frac{\varepsilon}{\log(n)})$ to drive the regularizing term $-\eta H(X)$ small , while the smoothness condition number of the above dual objective  is large at least at the exponential order of $\eta^{-1}$. On the other hand, the original UOT objective is non-smooth due to KL divergences. Therefore, direct application of gradient methods to solve either the primal or the dual of entropic regularized UOT would not result in competitive convergence rate, which perhaps is the reason for the limited literature of gradient methods in the context of UOT.
Recently, \citep{sparseOT_solomon, blondel2018smooth} also observe that entropic regularization imposed by the Sinkhorn algorithm  keeps the transportation plan dense and strictly positive. 
\end{remark}

\subsection{Squared $\ell_2$-norm Regularized UOT and Sparsity}

Remark \ref{entropic_downvote} motivates our choice of squared $\ell_2$-norm regularization, which would later result in better condition number in the dual formulation with only linear dependence on $\eta^{-1}$. Of practical interest, it also leads to sparse transportation plan empirically as later demonstrated in our experimental Section \ref{sparse_sec}.
\begin{remark}
While squared $\ell_2$-norm was well-studied as a sparsity-induced regularization in standard OT \citep{blondel2018smooth, peyre2020computational}, for which dedicated solvers have been developed \citep{lorenz2019quadratically, sparseOT_solomon}, to the best of our knowledge, no work has considered squared $\ell_2$-norm regularization in the context of UOT. This paper is the first to propose the squared $\ell_2$-norm  UOT problem, consequently derive feasible algorithms, and empirically investigate the sparsity of the output couplings. 
\end{remark}
The squared $\ell_2$-norm  UOT problem is formulated as: 
\begin{equation}
\label{sparse_UOT_func}
    \UOT^{\eta}_{\KL}(\aA, \bB) =  \min_{X \in \Br_{+}^{n \times n}} \{ g_\eta(X) := \left\langle C, X\right\rangle + \eta \| X\|_2^2 + \tau \KL(X \one_{n} || \aA) + \tau \KL(X^{\top} \one_{n} || \bB)\}, 
\end{equation}
where $\eta > 0$ is a given regularization parameter.
Let $ X^\eta  = \argminE_{X \in \Br_{+}^{n \times n}} g_\eta(X)  $, which is unique by the strong convexity of $g_\eta(X)$. 
We proceed to establish the duality of the squared $\ell_2$-norm UOT problem, which is based on the Fenchel-Rockafellar Theorem and would be the basis for our algorithmic development.



\begin{lemma}
\label{duality_theorem2}
The dual problem to \eqref{sparse_UOT_func} can be written as: 
\begin{align}
\nonumber
    \max_{\xX=(\uU, \vV) \in \Br^{2n}} \Big \{F_{a}(\xX) :=& -\frac{1}{4 \eta }\sum_{i,j=1}^{n} \max\{0, u_i+ v_j - C_{ij}\}^2 \\
    & - \tau \left\langle  e^{-\uU/\tau} , \aA\right\rangle- \tau \left\langle  e^{-\vV/\tau} , \bB\right\rangle +\tau \aA^{\top}  \one_{n} +\tau \bB^{\top}  \one_{n}  \Big\}.  \label{alt_dual_prob_adaptive1} 
\end{align}
Let $(\uU^*, \vV^*)$ be an optimal solution to \eqref{alt_dual_prob_adaptive1}, then the primal solution to \eqref{sparse_UOT_func} is given by:
\begin{align}
\label{alt_primal_sol}
    X^\eta_{ij} = \frac{1}{2\eta}\max\{0, u^*_i+ v^*_j - C_{ij}\}, \quad \forall i,j\in [n].
\end{align}
Moreover, we have  $\forall i,j\in[n]$:
\begin{align}
    -&\frac{u^*_i}{\tau} + \log(a_i) = \log(\sum_{k=1}^n X^\eta_{ik}), \label{alt_gradient_cond1}\\
    -&\frac{v^*_j}{\tau} + \log(b_j) = \log(\sum_{k=1}^n X^\eta_{kj}). \label{alt_gradient_cond2}
\end{align}
\end{lemma}
\begin{proof}
Given a convex set $S$, we consider the convex indicator function:
\begin{align*}
    \mathbbm{1}_S(X) = \begin{cases} 0 \quad \text{, if $X\in S$}\\ +\infty,  \text{otherwise}\end{cases},
\end{align*}
and write the primal objective \eqref{sparse_UOT_func} as:
\begin{equation}
   \UOT^{\eta}_{\KL}(\aA, \bB) =  \min_{X \in \Br^{n \times n}} \{ \left\langle C, X\right\rangle + \eta \| X\|_2^2 +\mathbbm{1}_S(X)+ \tau \KL(X \one_{n} || \aA) + \tau \KL(X^{\top} \one_{n} || \bB)\}  \label{dual_proof1}
\end{equation}
with $S = \{X\in \Br^{n \times n}| X_{ij}\geq 0, \forall i, j \in [n] \}$. Next, we consider the functions:
\begin{align*}
    G(X) = \frac{1}{\eta} \left\langle C, X \right\rangle  + \| X\|_2^2  +  \frac{1}{\eta} \Xx_S(X), \quad F_1(\yY) = \tau  \KL(\yY || \aA), \quad F_2(\yY) = \tau  \KL(\yY || \bB),
\end{align*}
whose convex conjugates are as follows:
\begin{align*}
    G^*(p)&:=\sup_{X\in \Br^{n\times n}}\{ \left\langle p, X \right\rangle - G(X)\} =\frac{1}{4} \sum_{i,j=1}^n \max\{0, p_{ij} -\frac{1}{\eta} C_{ij}  \}^2,\\
     F_1^*(\uU) &:= \sup_{\yY\in \Br^{n}}\{\left\langle \uU, \yY \right\rangle-  F_1(\yY) \}=  \tau \left\langle  e^{\uU/\tau} , \aA\right\rangle - \tau \aA^T  \one_{n},\\
    F_2^*(\vV) &:= \sup_{\yY\in \Br^{n}}\{\left\langle \vV, \yY \right\rangle-  F_2(\yY)\} =  \tau \left\langle  e^{\vV/\tau} , \bB\right\rangle - \tau \bB^T  \one_{n}.
\end{align*}
Consider the linear operator $A: \Br^{2n} \rightarrow \Br^{n\times n}$ that maps $A(\uU, \vV) = X$ with 
$X_{ij} = u_i + v_j$,
then $A$ is continuous and its adjoint $A^*: \Br^{n\times n} \rightarrow \Br^{2n}$ is $A^*(X) = (X \one_{n}, X^T \one_{n})$. 
Now note that the problem \eqref{alt_dual_prob_adaptive1} can be rewritten as:
\begin{align*}
     \max_{\uU, \vV \in \Br^{n}} - F_1^*(-\uU) - F_2^*(-\vV) - \eta G^*\bigg(\frac{A(\uU, \vV)}{\eta}\bigg).
\end{align*}
By the Fenchel-Rockafellar  Theorem (restated as Theorem \ref{thm:frenchel_dual} in Appendix \ref{appen_supple}), we obtain its   dual problem as: 
\begin{align*}
    \inf_{X \in \Br^{n \times n }} F_1(X \one_{n}) + F_2(X^T \one_{n}) + \eta G(X),
\end{align*}
which is the optimization problem \eqref{dual_proof1} and obtain that:
\begin{align*}
    &X^\eta \in \partial G^*\bigg(\frac{A(\uU^*, \vV^*)}{\eta}\bigg) \implies  \forall i, j \in [n]:  X^\eta_{ij} = \frac{1}{2\eta}\max\{0, u^*_i+ v^*_j - C_{ij}\}, \\ 
    &-\uU^* \in \partial F_1(X^\eta \one_{n}) \implies \forall i: -\frac{u^*_i}{\tau} + \log(a_i) = \log(\sum_{k=1}^n X^\eta_{ik}), \\
    &-\vV^* \in \partial F_2((X^\eta)^T \one_{n}) \implies \forall j: -\frac{v^*_j}{\tau} + \log(b_j) = \log(\sum_{k=1}^n X^\eta_{kj}),    
\end{align*}
which conclude the proof of the Lemma.
\end{proof}

An equivalent dual form of \eqref{sparse_UOT_func} can be deducted from the above Lemma \eqref{duality_theorem2} and presented in the next Corollary.


\begin{corollary}
\label{duality_theorem}
The dual problem to \eqref{sparse_UOT_func} can be written as: 
\begin{equation}
    \max_{(\uU, \vV, \tT) \in \Xx} \Big \{ -\frac{1}{4 \eta }\sum_{i,j=1}^{n} t_{ij}^2 -   \tau \left\langle  e^{-\uU/\tau} , \aA\right\rangle - \tau \left\langle  e^{-\vV/\tau} , \bB\right\rangle + \tau \aA^{\top}  \one_{n} + \tau \bB^{\top}  \one_{n} \Big \},  \label{dual_prob1} 
\end{equation}
where $\Xx = \{(\uU, \vV, \tT) | \uU, \vV \in \Br^{n}, \tT \in \Br^{n\times n}: t_{ij} \geq 0, t_{ij} \geq u_i + v_j - C_{ij} \quad \forall i, j \}$.
Let $(\uU^*, \vV^*, \tT^*)$ be an optimal solution to \eqref{dual_prob1}, then the primal solution to \eqref{sparse_UOT_func} is given by:
\begin{equation}
\label{primal_sol}
    X^\eta_{ij} = \frac{1}{2\eta}t^*_{ij}.
\end{equation}
Moreover, we have have $\forall i,j\in[n]$:
\begin{align}
    -&\frac{u^*_i}{\tau} + \log(a_i) = \log(\sum_{k=1}^n X^\eta_{ik}), \label{gradient_cond1}\\
     -&\frac{v^*_j}{\tau} + \log(b_j) = \log(\sum_{k=1}^n X^\eta_{kj}), \label{gradient_cond2}\\
    &t^*_{ij} = \max\{0, u^*_i+ v^*_j - C_{ij}\}.\label{gradient_cond3}
\end{align}
\end{corollary}
\begin{proof}
From the dual objective \eqref{alt_dual_prob_adaptive1} in Lemma   \ref{duality_theorem2},   we replace $ \max\{0,u_i + v_j - C_{ij}\}$ with the dummy variable $t_{ij}$ constrained by $t_{ij} \geq 0$ and $t_{ij}\geq u_i + v_j - C_{ij}$. We have just rewritten our  problem \eqref{sparse_UOT_func} as the optimization problem \eqref{dual_prob1} over the convex set $\Xx$. We thus have \eqref{gradient_cond3} at the optimal solution, while \eqref{primal_sol}, \eqref{gradient_cond1} and \eqref{gradient_cond2} directly follow from \eqref{alt_primal_sol}, \eqref{alt_gradient_cond1} and \eqref{alt_gradient_cond2}. 
\end{proof}

\begin{remark}
\label{equi_sol_remark}
We note that if $(\uU^*, \vV^*, \tT^*)$ is an optimal solution to \eqref{dual_prob1}, then $(\uU^*, \vV^*)$ is an optimal solution to \eqref{alt_dual_prob_adaptive1}.
\end{remark}



\section{Complexity Analysis of Approximating 
UOT}\label{sec_analysis}

In this Section, we provide the algorithmic development based on Gradient Extrapolation Method (GEM) for approximating UOT. 
We first define the quantities within interest  and present the regularity conditions in Section  \ref{sec_assumptions}.  Next, we prove several helpful properties of  squared $\ell_2$-norm UOT in Section \ref{sec:helper_uot}. We then characterize the problem as composite optimization in Section \ref{sec:dual}, on which GEM-UOT is derived.  The complexity analysis for  GEM-UOT  follows in Section \ref{sec:gem}.
In Section \ref{sec:gemruot}, we consider a relaxed UOT problem that is relevant to certain  practical applications, and present a easy-to-implement and practical algorithm, termed GEM-RUOT, for such setting.








\subsection{List of Quantities and Assumptions}
\label{sec_assumptions}
Given the two masses $\aA, \bB \in \Br^n_+$, we define the following notations and quantities:
\begin{align*}
&a_{min} = \min_{1 \leq i \leq n}\{a_i\}, \quad b_{min} = \min_{1 \leq j \leq n}\{b_j\}, \quad a_{max} = \|\aA\|_\infty, \quad b_{max} = \|\bB\|_\infty\\
 \kappa&= \frac{ 1}{\min\{a_{min}, b_{min}\}}, \quad R =\frac{(\alpha + \beta)^2}{4}, \quad p=\frac{1}{2} \min\{a_{min}, b_{min}\} e^{-\frac{D}{\tau}}, \quad q = \alpha + \beta, \\
\nonumber
 D&= \|C\|_\infty+   \eta( \alpha + \beta )  + \tau  \log\bigg( \frac{\alpha + \beta}{2} \bigg)  -\tau \min\{\log(a_{min}), \log(b_{min})  \},\\
\nonumber
L_1 &= \|C\|_\infty + 2 \eta q+  2\tau| \log(p)|+ 2\tau| \log(q)| + \tau \max_i |\log(a_i)|+ \tau \max_i |\log(b_i)|.
\end{align*}

We hereby present the assumptions (A1-A3) required by our algorithm. For interpretation, we also restate the regularity conditions of Sinkhorn \citep{pham2020unbalanced} for solving UOT, where detailed discussion on how their assumptions in the original paper are equivalent to (S1-S4) is deferred to Appendix \ref{sinkhorn_assumption}.

\begin{minipage}[t]{0.5\textwidth}
\textbf{Regularity Conditions of  this Paper}

(A1) $a_{min} >0$, $b_{min} > 0$.

(A2) $|\log(a_{min})| = O(\log(n)),$ \\ \text{  }\text{  }\text{  } $\quad |\log(b_{min})| = O(\log(n)) $.

(A3) $ \tau = \Omega(\min\{\frac{1}{\alpha + \beta}, \|C \|_\infty\})$. 

\end{minipage}
\hfill
\begin{minipage}[t]{0.5\textwidth}

\textbf{Regularity Conditions of Sinkhorn}

(S1) $a_{min} >0$, $b_{min} > 0$. 

(S2) $|\log(a_{min})| = O(\log(n)),$ \\ \text{  }\text{  }\text{  } $ \quad |\log(b_{min})| = O(\log(n)) $.

(S3)$|\log(a_{max})|=O(\log(n)),$ \\ \text{  }\text{  }\text{  } $ \quad  |\log(b_{max})| =O(\log(n))$.

(S4) $\alpha, \beta, \tau$ are positive constants.

\end{minipage}

\begin{remark}
\label{remark_assumption}
Compared to the regularity conditions of Sinkhorn, ours  lift the strict assumptions (S3) and (S4) that put an upper bound on $\tau$ and the input masses. Thus our complexity analysis supports much more flexibility for the input masses and is still suitable to applications requiring large $\tau$ to enforce the marginal constraints. On the other hand, our method requires very mild condition (A3) on $\tau$. For example, (A3) is naturally satisfied by (S4) combined with the boundedness of the grounded cost matrix $C$ which is widely assumed by the literature \citep{altschuler2017near, dvurechensky2018computational}. We strongly note that our method works for any $\alpha, \beta, \|C\|_\infty$ as long as $\tau$, a parameter under our control, satisfies (A3). 
\end{remark}











\subsection{Properties of Squared $\ell_2$-norm UOT}
\label{sec:helper_uot}
We next present several useful properties and related bounds for the squared $\ell_2$-norm UOT problem.

\begin{lemma}
\label{obj_identities}
The following identities hold: 
\begin{align}
g_\eta( X^\eta) + 2\tau \| X^\eta\|_1 + \eta\| X^\eta\|_2^2= \tau(\alpha+\beta), \label{obj_identities1} \\
f(X_f) + 2\tau \|X_f\|_1= \tau(\alpha+\beta). \label{obj_identities2}
\end{align}
\end{lemma}
\begin{proof}
The identity \eqref{obj_identities2}  follows from \citep[Lemma 4]{pham2020unbalanced}. 
To prove \eqref{obj_identities1}, we consider the function $g_\eta(t X^\eta)$, where $t \in \Br^{+}$, we have
\begin{align*}
   g_\eta(t X^\eta) =   \left\langle C, t X^\eta\right\rangle + \eta \| t X^\eta\|_2^2 + \tau \KL(t X^\eta \one_{n} || \aA) + \tau \KL((t X^\eta)^{\top} \one_{n} || \bB). 
\end{align*}
Simple algebraic manipulation gives:
\begin{align*}
    \KL(t X^\eta \one_{n} || \aA)&= t  \KL( X^\eta \one_{n} || \aA)+ (1-t) \alpha + \| X^\eta\|_1  t \log(t)\\
    \KL((t X^\eta)^{\top} \one_{n} || \bB) &= t  \KL(( X^\eta)^{\top} \one_{n} || \bB) + (1-t) \beta + \| X^\eta\|_1  t \log(t).
\end{align*}
We thus obtain that:
\begin{align*}
     g_\eta(t X^\eta) = t g_\eta( X^\eta) +\tau(1-t)(\alpha + \beta) + 2\tau \| X^\eta\|_1  t \log(t) + \eta(t^2-t)\| X^\eta\|_2^2.
\end{align*}
Differentiating $g_\eta(t X^\eta)$ with respect to $t$:
\begin{align*}
    \frac{\partial  g_\eta(t X^\eta)}{\partial t} = g_\eta( X^\eta) - \tau(\alpha+\beta) + 2\tau \| X^\eta\|_1  (1+ \log(t)) +  \eta(2t-1)\| X^\eta\|_2^2.
\end{align*}
From the above analysis, we can see that $g_\eta(t X^\eta)$ is well-defined for all $t \in \Br^{+}$ and attains its minimum at $t=1$. Setting $ \frac{\partial  g_\eta(t X^\eta)}{\partial t}\at{t=1}=0$, we obtain the identity \eqref{obj_identities1}.
\end{proof}
The next Lemma provides bounds for both the primal and dual solutions to the squared $\ell_2$-norm UOT problem.

\begin{lemma}
\label{opt_sol_bounds}
We have the following bounds for the optimal solution $(\uU^*, \vV^*, \tT^*)$ of \eqref{dual_prob1}:
\begin{align}
\| X^\eta\|_1 &\leq \frac{\alpha + \beta}{2},   \label{opt_sol_bounds1} \\
\|X_f\|_1 &\leq \frac{\alpha + \beta}{2},    \label{opt_sol_bounds2}\\
 u_i^* &\geq   \tau \log \left(\frac{2a_i}{\alpha+\beta}\right) \forall i \in [n],\\
 v_j^* &\geq   \tau \log \left(\frac{2b_j}{\alpha+\beta}\right) \forall j \in [n], \\
\|\uU^*\|_\infty, \|\vV^*\|_\infty &\leq D, \label{opt_sol_bounds3}\\
\min_{i,j} \Big \{ \sum_{k=1}^n X^\eta_{ik}, \sum_{k=1}^n X^\eta_{kj} \Big \} &\geq \min\{a_{min}, b_{min}\} e^{-\frac{D}{\tau}}. \label{opt_sol_bounds4}
\end{align}
\end{lemma}
\begin{proof}
Noting that $f(X)$ and $g_\eta(X)$ are non-negative for $X\in\Br_{+}^{n \times n}$, we directly obtain \eqref{opt_sol_bounds1} and \eqref{opt_sol_bounds2} from Lemma \ref{obj_identities}.
By $ X^\eta_{ij}\geq0$, $\forall i, j = 1,\dots,n$, and  \eqref{opt_sol_bounds1}, we have:
\begin{align}
    \frac{\alpha + \beta}{2} \geq \| X^\eta\|_1 =  \sum_{i,j=1}^n X^\eta_{ij} \geq  \max\{\sum_{k=1}^n X^\eta_{ik},  \sum_{k=1}^n X^\eta_{kj} \}. \label{eq_02}
\end{align}
We obtain in view of Lemma \ref{duality_theorem2}  that: 
\begin{align}
        \log (a_i) - \frac{u_i^*}{\tau} = \log \big(  \sum_{j=1}^n X^\eta_{ij} \big) \leq \log\big( \frac{\alpha + \beta}{2} \big) , \ i = 1,\dots,n , \label{eq_01a} \\
       \log (b_j) - \frac{v_j^*}{\tau} = \log \big(  \sum_{i=1}^n X^\eta_{ij} \big) \leq   \log\big( \frac{\alpha + \beta}{2} \big), \ j = 1,\dots,n, \label{eq_01b}
\end{align}
which are equivalent to:
\begin{align}
    u_i^* \geq   \tau \log \left(\frac{2a_i}{\alpha+\beta}\right) \geq \tau \left( \log (a_{min}) -  \log\big( \frac{\alpha + \beta}{2} \big) \right) , \ i = 1,\dots,n,  \label{eq_03a}  \\
    v_j^* \geq   \tau \log \left(\frac{2b_j}{\alpha+\beta}\right) \geq \tau \left ( \log (b_{min}) -  \log\bigg( \frac{\alpha + \beta}{2} \bigg) \right) , \ j = 1,\dots,n.  \label{eq_03b}
\end{align}
By \eqref{primal_sol}, 
we have $\frac{1}{2\eta}(u_i^* + v_j^*- C_{ij}) \leq  X^\eta_{ij} $, $\forall i,j = 1,\dots,n$. Hence, 
\begin{align}
    u_i^* \leq  2 \eta X^\eta_{ij} + C_{ij} - v_j^* \overset{\eqref{eq_03b}}{\leq} 2 \eta X^\eta_{ij} + C_{ij} - \tau \left ( \log (b_{min}) -  \log\big( \frac{\alpha + \beta}{2} \big) \right),  \label{eq_04a} \\
    v_j^* \leq 2 \eta X^\eta_{ij} + C_{ij} - u_i^* \overset{\eqref{eq_03a}}{\leq} 2 \eta X^\eta_{ij} + C_{ij} - \tau \left( \log (a_{min}) - \log\big( \frac{\alpha + \beta}{2} \big) \right).  \label{eq_04b}
\end{align}
Note that $X_{ij}^\eta \leq \max_{i,j} | X_{ij}^\eta| = \| X^\eta \|_\infty \leq \| X^\eta \|_1 \leq \frac{\alpha + \beta}{2}$. We have
\begin{align*}
   \tau \left( \log (a_{min}) -  \log\bigg( \frac{\alpha + \beta}{2} \bigg) \right)&\leq u_i^* \leq \eta( \alpha + \beta )+ C_{ij} - \tau \left ( \log (b_{min}) -  \log\bigg( \frac{\alpha + \beta}{2} \bigg) \right), \\
   \tau \left( \log (b_{min}) -  \log\bigg( \frac{\alpha + \beta}{2} \bigg) \right) &\leq v_j^* \leq  \eta( \alpha + \beta ) + C_{ij} - \tau \left( \log (a_{min}) - \log\bigg( \frac{\alpha + \beta}{2} \bigg) \right),\\
\therefore \max\{\|\uU^*\|_\infty , \|\vV^*\|_\infty \} &\leq D = \|C\|_\infty+   \eta( \alpha + \beta )  + \tau  \log\big( \frac{\alpha + \beta}{2} \big)  \\
&\quad -\tau \min\{\log(a_{min}), \log(b_{min}) \}.
\end{align*}
To prove \eqref{opt_sol_bounds4}, we note that by Lemma \ref{duality_theorem2}, $\forall i: \sum_{k=1}^{n} X^\eta_{ik} =  a_i \cdot e^{\frac{-u_i^*}{\tau}}\geq  a_i \cdot  e^{\frac{-\|\uU^*\|_\infty}{\tau}} \geq  a_i \cdot e^{\frac{-D}{\tau}}$.
Similarly, we have $\forall j: \sum_{k=1}^{n} X^\eta_{kj} \geq b_j \cdot  e^{\frac{-D}{\tau}}$ and thus obtain that:
$$\min_{i,j}\{ \sum_{k=1}^n X^\eta_{ik}, \sum_{k=1}^n X^\eta_{kj} \} \geq \min\{a_{min}, b_{min}\} e^{-\frac{D}{\tau}},$$
which concludes the proof of this Lemma.
\end{proof}

The local Lipschitzness of the squared $\ell_2$-norm UOT objective is established in the following Lemma.

\begin{lemma}
\label{sparseUOT_lipschitz}
For $0<p\leq q$ defined in Section \ref{sec_assumptions}, let $U_{p,q} = \{ X \in \Br_{+}^{n \times n}| \forall i, j: \sum_{k=1}^n X_{ik} \in [p, q], \sum_{k=1}^n X_{kj} \in [p, q]  \} $. 
Then $g_\eta(X)$ is $L_g$-Lipschitz in $U_{p,q}$, i.e.  $\forall X,Y \in U_{p,q}$:
\begin{align*}
    |g_\eta(X) - g_\eta(Y)| \leq L_1 \|X-Y\|_1,
\end{align*}
where $L_1 = \|C\|_\infty + 2\eta q + 2\tau |\log(p)|+2\tau |\log(q)| + \tau \max_i |\log(a_i)|+ \tau \max_i |\log(b_i)| $.
\end{lemma}
\begin{proof}
Take any $X, Y \in U_{p,q}$. By Mean Value Theorem, there exists $c\in (0,1)$ and $Z= (1-c)X + c Y \in U_{p,q}$ such that: 
\begin{align}
    |g_\eta(X) - g_\eta(Y)| &= |\langle \nabla g_\eta(Z),  X - Y\rangle| \leq \|\nabla g_\eta(Z) \|_\infty \|X-Y\|_1 , \label{lipschitz_inter1}
\end{align}
where we have used Holder inequality.
We also have:
\begin{align*}
    \frac{\partial g_\eta}{\partial Z_{ij}} = C_{ij} + 2 \eta Z_{ij} + \tau \log(\sum_{k=1}^n Z_{ik}) + \tau \log(\sum_{k=1}^n Z_{kj}) -\tau \log(a_i) -\tau \log(b_j).
\end{align*}
Since $Z \in U_{p,q}$, we obtain $\forall i,j: \sum_{k=1}^n Z_{ik}, \sum_{k=1}^n Z_{kj} \in [p, q]$, which implies that:
\begin{align*}
    &|\log(\sum_{k=1}^n Z_{ik})|, | \log(\sum_{k=1}^n Z_{kj})| \leq | \log(p)| + |\log(q)|
\end{align*}
From the above observations, we consequently obtain that:
\begin{align}
\nonumber 
    \|\nabla g_\eta(Z) \|_\infty  &\leq \|C\|_\infty + 2\eta q +2\tau |\log(p)|+2\tau |\log(q)| \\
    &\quad + \tau \max_i |\log(a_i)|+ \tau \max_i |\log(b_i)| = L_g \label{lipschitz_inter2}
\end{align}
Plugging \eqref{lipschitz_inter2} into \eqref{lipschitz_inter1} yields the desired inequality.
\end{proof}

\subsection{Characterization of the Dual Objective}\label{sec:dual}
GEM-UOT is designed by solving the dual  \eqref{dual_prob1} from Corollary \ref{duality_theorem}, which is equivalent to:
\begin{equation}
\label{dual_sparse_UOT}
     \min_{\xX \in \Xx} h_\eta(\xX = (\uU, \vV, \tT)) := \frac{1}{4 \eta }\sum_{i,j=1}^{n} t_{ij}^2 +   \tau \left\langle  e^{-\uU/\tau} , \aA\right\rangle + \tau \left\langle  e^{-\vV/\tau} , \bB\right\rangle. 
\end{equation}
For $\xX = (\uU, \vV, \tT)$, we consider the functions: 
\begin{align}
    f_\eta(\xX &) = \tau \left\langle  e^{-\uU/\tau} , \aA\right\rangle + \tau \left\langle  e^{-\vV/\tau} , \bB\right\rangle -\frac{\min\{a_{min}, b_{min}\}}{2\tau}  e^{-D/\tau} \big(\|\uU\|_2^2+\|\vV\|_2^2\big) , \label{f2}\\
    w_\eta(\xX&) = \frac{\min\{a_{min}, b_{min}\}}{2\tau}  e^{-D/\tau} \big(\|\uU\|_2^2+\|\vV\|_2^2\big) +\frac{1}{4 \eta }\sum_{i,j=1}^{n} t_{ij}^2.     \label{w_eta}
\end{align}
Then the problem \eqref{dual_sparse_UOT} can be rewritten as:
 \begin{align}
     \min_{\xX \in \Xx} h_\eta(\xX = (\uU, \vV, \tT)) := f_\eta(\xX) + w_\eta(\xX), \label{primal-dual-form}
 \end{align}
which can be characterized as the composite optimization over the sum of the locally convex smooth $f_\eta(\xX)$ and the locally strongly convex $w_\eta(\xX)$ by the following Lemma \ref{sparse_UOT_lipschitz_lemma}, whose proof is deferred to Appendix \ref{lma10proof}. 

\begin{lemma}
\label{sparse_UOT_lipschitz_lemma}
Let  $V_D = \{ \xX  = (\uU, \vV, \tT) \in \Br^{n^2+ 2n} :  \forall i, j \in [n],  \tau \log(\frac{2a_i}{\alpha+\beta}) \leq u_i \leq D,  \tau \log(\frac{2b_j}{\alpha+\beta}) \leq  v_j \leq D\}$. Then $f_\eta(\xX)$ is  convex and $L$-smooth in the domain $V_D$, and $w_\eta(\xX)$ is $\mu$-strongly convex with:
\begin{align*}
    L &= \frac{\alpha+\beta}{2\tau}+\frac{\min\{a_{min}, b_{min}\}}{\tau}  e^{-D/\tau}, \\
     \mu &=  \min\bigg\{  \frac{\min\{a_{min}, b_{min}\}}{\tau}  e^{-D/\tau}, \frac{1}{2\eta} \bigg\}.
\end{align*}
\end{lemma}

The choice of $D$ in the above characterization is motivated by the next Corollary \ref{dual_domain}, which ensures that the optimal solution $\xX^*$ lies within the locality $V_D$. 
\begin{corollary}
\label{dual_domain}
If $\xX^*$ is an optimal solution to \eqref{primal-dual-form}, then $\xX^* \in V_D$.
\end{corollary}
\begin{proof}
The result follows directly from Lemma \ref{opt_sol_bounds}.
\end{proof}


\subsection{Gradient Extrapolation Method for UOT (GEM-UOT)}\label{sec:gem}

\textbf{Algorithm Description:}
We develop our algorithm based on Gradient Extrapolation Method (GEM), called GEM-UOT, to solve the  UOT problem. GEM was proposed by \citep{Lan2018RandomGE} to address problem in the form of \eqref{primal-dual-form} and viewed as  the dual of the Nesterov’s accelerated gradient  method. 
We hereby denote the prox-function associated with $w_\eta$ by:
\begin{align*}
P(\xX_0,\xX) := \frac{w_\eta(\xX) - [w_\eta(\xX_0) + \langle \nabla w_\eta(\xX_0), \xX-\xX_0 \rangle ]}{\mu},
\end{align*}
and the  proximal mapping associated with the closed convex set $\Yy$ and $w$ is defined as:
\begin{align}
    \mathcal{M}_\Yy(\gG, \xX_0, \theta) := \argminE_{\xX = (\uU,\vV, \tT) \in \Yy} \{\langle \gG,\xX \rangle + w_\eta(\xX) + \theta P(\xX_0,\xX) \}. \label{ulti_dual_problem}
\end{align}
One challenge for optimizing \eqref{primal-dual-form} is the fact that the component function $f_\eta(\xX)$ is smooth and convex  only on the locality $V_D$ as in Lemma \ref{sparse_UOT_lipschitz_lemma}. By Corollary \ref{dual_domain}, we can further write \eqref{primal-dual-form} as:
\begin{align}
\min_{\xX \in \Xx} h_\eta(\xX = (\uU, \vV, \tT)) =  \min_{\xX \in V_D \cap \Xx } h_\eta(\xX = (\uU, \vV, \tT)).  \label{equi_dual_obj1}
\end{align}
While minimizing the UOT dual objective, GEM-UOT hence projects its iterates onto the closed convex set $V_D \cap \Xx $ through the proximal mapping $\mathcal{M}_{V_D \cap \Xx}()$ to enforce the smoothness and convexity of  $f_\eta(\xX)$ on the convex hull of all the iterates, while ensuring they are in the feasible domain. Finally, GEM-UOT  outputs the approximate solution,  computed from the iterates, to the original UOT problem. The steps are summarized in Algorithm \ref{alg: acc-proximal-gd-uot}.

\begin{minipage}[t]{0.46\textwidth}
\begin{algorithm}[H]  
\caption{GEM-UOT}
\label{alg: acc-proximal-gd-uot}
\begin{algorithmic}[1] 
\STATE {\bfseries Input:} $ C, \aA, \bB, \varepsilon, \tau, \eta$
\STATE {\bfseries Initialization:}  $\underline{\xX}^0 = \xX^0 = \mathbf{0}, \quad \yY^{-1}=\yY^{0}= \mathbf{0}$ 
\STATE Compute $L, \mu$ based Lemma \ref{sparse_UOT_lipschitz_lemma}
  \STATE $\zeta = 1- 1/\big(1+ \sqrt{1+16 L/\mu}\big)$
  \STATE $\psi = \frac{1}{1-\zeta} -1, \quad \rho = \frac{\zeta}{1-\zeta}\mu$
\FOR{$t = 0, 1, 2, ..., k$}
\STATE $\widetilde{\yY}^t = \yY^{t-1} + \zeta (\yY^{t-1} - \yY^{t-2})$ \alglinelabel{iterate1}
\STATE  $\xX^t = \mathcal{M}_{V_D\cap \Xx}(\widetilde{\yY}^t, \xX^{t-1}, \rho)$ \alglinelabel{projection_step1}
\STATE  $\underline{\xX}^t = (1+\psi)^{-1} (\xX^t + \psi \underline{\xX}^{t-1} )$ \alglinelabel{iterate2}
\STATE  $\yY^t = \nabla f_\eta(\underline{\xX}^t)$ \alglinelabel{GEM_gradient}
\STATE  $\theta_t = \theta_{t-1}\zeta^{-1}$
\ENDFOR
\STATE Compute $\underline{\xX}^k = (\sum_{t=1}^k \theta_t )^{-1} \sum_{t=1}^k \theta_t \xX^t$ 
\STATE  Let $\underline{\uU}^k, \underline{\vV}^k ,\underline{\tT}^k$ be the vectors such that $\underline{\xX}^k = (\underline{\uU}^k, \underline{\vV}^k, \underline{\tT}^k) $
\FOR{$i,j = 1 \rightarrow n$}
\STATE $X^k_{ij}= \frac{1}{2\eta} \underline{t}^k_{ij}$ \alglinelabel{convert_dual_primal} 
\ENDFOR
\STATE {\bfseries Return} $X^k$
\end{algorithmic}
\end{algorithm}
\end{minipage}
\hfill
\begin{minipage}[t]{0.46\textwidth}
\begin{algorithm}[H]
\caption{GEM-RUOT}
\label{alg: acc-proximal-gd-uot-convex}
\begin{algorithmic}[1] 
\STATE {\bfseries Input:} $ C, \aA, \bB, \varepsilon, \tau, \eta$
\STATE {\bfseries Initialization:}   $\underline{\xX}^0 = \xX^0 = \mathbf{0}, \quad \yY^{-1}=\yY^{0}= \mathbf{0}$ 
\STATE Compute $L_a$ based on Lemma \ref{rgem_convex_lipschitz_lemma}
  
\FOR{$t = 0, 1, 2, ..., k$}
\STATE $\zeta_t = \frac{t-1}{t}, \quad \psi_t = \frac{t-1}{2}, \quad \rho_t = \frac{6L_a}{t},$
\STATE $\widetilde{\yY}^t = \yY^{t-1} + \zeta_t (\yY^{t-1} - \yY^{t-2})$ \alglinelabel{conv_iterate1}
\STATE  $\xX^t = \bar{\mathcal{M}}_{V_a}(\widetilde{\yY}^t, \xX^{t-1}, \rho_t)$ \alglinelabel{conv_projection_step1}
\STATE  $\underline{\xX}^t = (1+\psi_t)^{-1} (\xX^t + \psi_t \underline{\xX}^{t-1} )$ \alglinelabel{conv_iterate2}
\STATE  $\yY^t = \nabla h_a(\underline{\xX}^t)$ \alglinelabel{conv_GEM_gradient}
\STATE  $\theta_t = \theta_{t-1}\zeta_t^{-1}$
\ENDFOR
\STATE  $ \underline{\xX}^k = (\sum_{t=1}^k \theta_t )^{-1} \sum_{t=1}^k \theta_t \xX^t$ 
\STATE  Let $\underline{\uU}^k, \underline{\vV}^k$ be the vectors such that $\underline{\xX}^k = (\underline{\uU}^k, \underline{\vV}^k) $
\FOR{$i,j = 1 \rightarrow n$}
\STATE $X^k_{ij}=  \frac{1}{2\eta} \max\{0, \underline{u}^k_{i} + \underline{v}^k_{j} - C_{ij}\}$ 
\ENDFOR
\STATE Compute $F_a(\underline{\xX}^k)$
\STATE {\bfseries Return} $F_a(\underline{\xX}^k), X^k$
\end{algorithmic}
\end{algorithm}
\end{minipage}
\text{}\\

\textbf{Complexity Analysis:} We proceed to establish the complexity for approximating UOT. The following quantities are used in our  analysis: 
\begin{align*}
\sigma_0^2 = \| \nabla f_\eta(\xX^0)\|_2^2, \quad \Delta_{0, \sigma_0} = \mu P(\xX^0, \xX^*) + h_\eta(\xX^0) -h_\eta(\xX^*) + \frac{\sigma_0^2}{\mu}. 
\end{align*}

\noindent
In the next Theorem, we further quantify the number of iterations
\begin{align*}
    K_0 &=  4 \left( 1+\sqrt{1+16L/\mu} \right)
     \log\bigg(\tfrac{4n^{2} \Delta_{0, \sigma_0}^{1/2} }{\eta  \mu } \max\{\tfrac{L_1}{\varepsilon}, \tfrac{e^{D/\tau}}{\min\{a_{min}, b_{min}\}}, \tfrac{1}{\alpha+\beta} \}\bigg)
\end{align*}
that is required for GEM-UOT to output the $\varepsilon$-approximate solution to $\UOT_\KL(\aA, \bB)$.


\begin{theorem}
\label{thm: acce_proximal_guarantee1}
For $\eta = \frac{\varepsilon}{2R}$ and $k \geq K_0$,  the output $X^k$ of Algorithm~\ref{alg: acc-proximal-gd-uot}  is the $\varepsilon$-approximate solution to the UOT problem, i.e.
\begin{align}
 f(X^k)  - f(X_f) = f(X^k) - \UOT_{\KL}(\aA, \bB) \leq \varepsilon.   \label{eq:maintheorem}  
\end{align}
\end{theorem}
\begin{proof}
Recalling that $ X^\eta  = \argminE_{X \in \Br_{+}^{n \times n}} g_\eta(X)  $, we get:
\begin{align}
\nonumber
    f(X^k)  - f(X_f) &= g_\eta(X^k)  - \eta \| X^k \|^2_2 - g_\eta(X_f)  + \eta \| X_f \|^2_2 \\
    &\leq (g_\eta(X^k) - g_\eta( X^\eta) +\eta (\| X_f \|^2_2 - \| X^k \|^2_2  ). \label{main_thrm_proof1}
\end{align}
From  Lemma \ref{opt_sol_bounds}, we have  $\|X_f\|_2^2 \leq \|X_f\|_1^2\leq \frac{(\alpha + \beta)^2}{4} = R$. Since $\eta = \frac{\varepsilon}{2R}$, we further obtain that: 
\begin{align}
\eta (\| X_f \|^2_2 - \| X^k \|^2_2) \leq \frac{\varepsilon}{2}. \label{main_thrm_proof2}
\end{align}
Plugging \eqref{main_thrm_proof2} into \eqref{main_thrm_proof1}, we have:
\begin{align}
      f(X^k)  - f(X_f) &\leq g_\eta(X^k) - g_\eta( X^\eta)  + \frac{\varepsilon}{2} . \label{main_thrm_proof3}
\end{align}
Thus, it is left to bound $g_\eta(X^k) - g_\eta( X^\eta)$. Next, we proceed establish the bound on $\| X^k - X^\eta\|_1$ in the following Lemma \ref{lma: acce_proximal_guarantee1}.

\begin{lemma}
\label{lma: acce_proximal_guarantee1}
If $k \geq K_0$, the output $X^k$ of Algorithm \ref{alg: acc-proximal-gd-uot} satisfies:
\begin{align}
\|X^k - X^\eta\|_1 \leq \min\bigg\{\frac{\varepsilon}{2 L_1} , \frac{\min\{a_{min}, b_{min}\} e^{-D/\tau}}{2}, \frac{\alpha+\beta}{2} \bigg\},  \label{eq: gem_lemma1.1}
\end{align}
\end{lemma}
The proof of Lemma \ref{lma: acce_proximal_guarantee1} can be found in Appendix \ref{appen_acce_proximal_guarantee1} where we first bound the quantity in terms of the distance $ \| \underline{\xX}^k - \xX^* \|_2$ between the dual iterate and dual solution, which provably vanishes to the desired error $\frac{\varepsilon}{2L_1}$ for large enough $k$ thanks to the guarantee from GEM being applied to the dual objective  \eqref{ulti_dual_problem}. Lemma \ref{lma: acce_proximal_guarantee1} also helps in deriving the next Lemma \ref{baby_lemma1}, whose proof is in Appendix \ref{appen_baby_lemma1}. 
\begin{lemma}
\label{baby_lemma1}
For $k\geq K_0$, we have $X^k, X^\eta \in U_{p,q}$.
\end{lemma}
Lemma \ref{baby_lemma1} is necessary  to invoke Lemma \ref{sparseUOT_lipschitz} on the Lipschitzness of $g_\eta(X)$ within the convex domain $U_{p,q}$. In particular, we obtain from Lemma \ref{sparseUOT_lipschitz} with $X^k, X^\eta \in U_{p,q}$ that:
\begin{align}
    g_\eta(X^k) - g_\eta( X^\eta) &\leq L_1 \| X^k - X^\eta\|_1 \label{main_thrm_proof4}
    \end{align}
From \eqref{eq: gem_lemma1.1} and \eqref{main_thrm_proof4}, we get:

\begin{align}
     g_\eta(X^k) - g_\eta( X^\eta) &\leq \frac{\varepsilon}{2}. \label{main_thrm_proof5}
\end{align}
By plugging \eqref{main_thrm_proof5} into  \eqref{main_thrm_proof3}, we conclude the desired inequality \eqref{eq:maintheorem}. 
\end{proof}

Finally, the  complexity of GEM-UOT is given by the following Corollary \ref{acceGD_complexity1}.
\begin{corollary}
\label{acceGD_complexity1}
Under the  Assumptions (A1-A3), Algorithm \ref{alg: acc-proximal-gd-uot} has the iteration complexity of 
${O}\left(  (\alpha+\beta){\kappa}   \log\left(\frac{\tau \cdot n (\alpha+\beta) }{\varepsilon}\right) \right)$ and $\widetilde{O}(n^2)$ cost per iteration\footnote{$\widetilde{O}()$ excludes the logarithmic terms in the complexity.}.
\end{corollary}
\begin{proof}
Note that Algorithm \ref{alg: acc-proximal-gd-uot} sets $\eta = \frac{\varepsilon}{2R} = \frac{2\varepsilon}{(\alpha+\beta)^2}$. 
Under the Assumptions (A1-A3), we have:
\begin{align*}
\frac{D}{\tau} &= \frac{\|C\|_\infty}{\tau}+  \frac{\varepsilon}{\tau( \alpha + \beta) } + \log\bigg( \frac{\alpha + \beta}{2} \bigg)    -\min\{\log(a_{min}), \log(b_{min})\} \\
&= O(1) + \log\bigg( \frac{\alpha + \beta}{2} \bigg)   - \min\{\log(a_{min}), \log(b_{min})\}
\end{align*}
By Lemma \ref{sparse_UOT_lipschitz_lemma} and under the Assumptions (A1-A2), we obtain that: 
\begin{align*}
    L &= O\bigg(\frac{\alpha+\beta}{2\tau}\bigg)\\
    \mu&= \Omega\bigg(\min\bigg\{\frac{(\alpha+\beta)^2}{4\varepsilon}, \frac{1}{\tau} \min\{a_{min}, b_{min}\} e^{-\log\big( \frac{\alpha + \beta}{2} \big) +\min\{\log(a_{min}), \log(b_{min})\}} \bigg\}\bigg)\\
    &=  \Omega \bigg( \min\bigg\{ \frac{(\alpha+\beta)^2}{4\varepsilon}, \frac{1}{\tau} \min\{a_{min}, b_{min}\}^2 \bigg\}\bigg)  \\
     &\geq  \Omega \bigg( \min\bigg\{ \frac{n^2 ( 2\min\{a_{min}, b_{min}\} )^2}{4\varepsilon},  \frac{2\min\{a_{min}, b_{min}\}^2}{\tau(\alpha+\beta)} \bigg\}\bigg)  \\
    &=  \Omega \bigg(  \frac{\min\{a_{min}, b_{min}\}^2}{\tau (\alpha+\beta)}\bigg)  \text{ (since $\tau = \Omega(\frac{1}{\alpha+\beta})$)}\\
     \sqrt{\frac{L}{\mu}} &= O\bigg( \frac{ (\alpha +\beta)}{\min\{a_{min}, b_{min}\}}\bigg) = O( (\alpha +\beta) {\kappa}).
\end{align*}
For convenience, we recall that $P(\xX_0,\xX) = \big(w_\eta(\xX) - [w_\eta(\xX_0) + \langle \nabla w_\eta(\xX_0), \xX-\xX_0 \rangle ]\big)/\mu$. 
From the initialization $\xX^0 = (\uU^0, \vV^0, \tT^0) = \mathbf{0}$, $\eta = \frac{\varepsilon}{2R} = \frac{2\varepsilon}{(\alpha+\beta)^2}$  and Lemma \ref{opt_sol_bounds}, we have:
\begin{align*}
    \frac{D}{\tau} &= O(1)  + \log\left( \frac{\alpha + \beta}{2} \right)- \min\{\log(a_{min}), \log(b_{min})\} \leq O(\log(n(\alpha+\beta))) \\
    \sigma_0^2 &= \| \nabla f_\eta(\xX^0)\|_2^2 = \|\aA\|_2^2 + \|\bB\|_2^2 \leq (\alpha + \beta)^2\\
    h_\eta(\xX^0) &= \tau (\alpha + \beta) \geq h_\eta(\xX^*) \geq 0,\quad w_\eta(\xX^0) = 0,\quad  \nabla w_\eta(\xX^0) = \mathbf{0}\\
    0 \leq w_\eta(\xX^*) &\leq  \frac{\min\{a_{min}, b_{min}\}}{2\tau}  \big(n\|\uU^*\|_\infty^2+n\|\vV\|_\infty^2\big) +\frac{1}{4 \eta } \|X^\eta\|_2^2 \\
    &\leq O(\min\{a_{min}, b_{min}\} \tau n  \log(n)^2 + (\alpha+\beta)^4/\varepsilon) \\
    \Delta_{0, \sigma_0} &= \mu P(\xX^0, \xX^*) + h_\eta(\xX^0) -h_\eta(\xX^*) + \frac{\sigma_0^2}{\mu} \\
    &= O\bigg({\tau \min\{a_{min}, b_{min}\} n\log(n)^2 } + \tau(\alpha+\beta) + \frac{\tau (\alpha+\beta)^2}{\min\{a_{min}, b_{min}\}^2}\bigg).
\end{align*}
We thus obtain the  asymptotic complexity under (A1-A3) for the following term inside $K_0$:
\begin{align*}
    \log\bigg(\tfrac{4n^{2} \Delta_{0, \sigma_0}^{1/2} }{\eta \mu } \max\{\tfrac{L_1}{\varepsilon}, \tfrac{e^{D/\tau}}{\min\{a_{min}, b_{min}\}}, \tfrac{1}{\alpha+\beta} \}\bigg) = O\bigg( \log\left(\frac{\tau \cdot n (\alpha+\beta)  }{\varepsilon}\right)\bigg).
\end{align*}
 Combining with the bounds of $L$ and $\mu$ above, we obtain that the iteration complexity is $K_0 = O\bigg( (\alpha +\beta) \kappa \log\left(\frac{\tau \cdot n (\alpha+\beta)  }{\varepsilon}\right)\bigg) $ under the assumptions (A1-A3). 

Lines \ref{iterate1} and \ref{iterate2} are  $O(1)$ point-wise update steps of vectors with $n^2 + 2n$ entries, thereby incurring $O(n^2)$ complexity. 
The proximal mapping on line \ref{projection_step1} is equivalent to solving a sparse separable quadratic program of size $O(n^2)$ with sparse linear constraints and  can be solved efficiently in $\widetilde{O}(n^2)$ time and with $O(n^2)$ space (see Appendix \ref{appen_projection_gemuot} for details). 
The computation of the gradient on line \ref{GEM_gradient} takes $O(n^2)$ since the partial derivative  $\partial f_\eta$ with respect to any of its $n^2 +2n$ variables takes $O(1)$ operations. Thus, the cost per iteration is $\widetilde{O}(n^2)$ in total.
\end{proof}

The complexities of GEM-UOT and other algorithms in the literature are summarized in Table \ref{complexity_table}.
GEM-UOT is asymptotically better than Sinkhorn in $\tau$ and $\varepsilon$, and, as mentioned before, has the edge on practicality over Sinkhorn. 
\begin{table*}[h]
\caption{Complexity of algorithms for solving UOT problems.}
\label{complexity_table}
  \centering
  {\renewcommand{\arraystretch}{1.6}
  \begin{tabular}{|p{2.7cm}|p{2.8cm}| p{1.7cm} |p{6.0cm}|}
    \hline 
    \textbf{Algorithms} &   \textbf{Assumptions}  &  \textbf{Cost per iteration} & \textbf{Iteration complexity} \\
    \hline
     MM  \citep{Chapel2021UnbalancedOT} &  Unspecified iteration complexity   &  ${O}\big(n^{2} \big) $& ${O}\big(n^{1.27} \big) $- Empirical Complexity \\
    \hline 
   
    Sinkhorn   \citep{pham2020unbalanced} & (S1), (S2), (S3), (S4) & 
    $O(n^2)$& ${O}\Big( { (\alpha + \beta) }  \frac{\tau \log(n)}{\varepsilon}   \log\left(\frac{ \log(n)(\alpha+\beta)}{\varepsilon} \right) \Big)$    \\
    \hline 
  GEM-UOT (this paper)  &   (A1)$\equiv$(S1), (A2)$\equiv$(S2), (A3) & $\widetilde{O}(n^2)$ &  $O\Big(  (\alpha + \beta){\kappa} \log\left(\frac{\tau \cdot n (\alpha+\beta)  }{\varepsilon}\right)  \Big) $   \\
    \hline
  \end{tabular}
  }
\end{table*}

   

\textbf{Novelty and Practicality:} GEM-UOT is predicated on the novel dual formulation of the  squared $\ell_2$-norm regularized UOT
and the intricate design of  smooth locality $V_D$ (Lemma \ref{sparse_UOT_lipschitz_lemma}). We note that naive application of GEM would not lead to competitive complexity. Beside function decomposition, we design the smooth locality $V_D$  to enforce tight bound for the condition number $L$. The techniques can benefit follow-up work tackling squared $\ell_2$-norm UOT. Next, we highlight the favorable features of GEM-UOT that Sinkhorn is short of: 
\begin{itemize}
    \item \textit{Compatibility with modern applications:} Gradient methods are specifically congruent with and thus deployed as heuristics without theoretical guarantee in many  emerging UOT applications   \citep{balaji2020robust, Yang_Scalable_2019, wassloss}. Applications \citep{PITIE2007123, courty2016optimal, muzellec2016tsallis} in which only sparse transportation plan  is of interest also nullify the usage of Sinkhorn, since entropic regularization imposed by it keeps the transportation plan dense and strictly positive \citep{sparseOT_solomon, blondel2018smooth}. GEM-UOT addresses both of these problems by being the first gradient-based sparse UOT solver in the literature.
    \item \textit{Logarithmic dependence on $\tau$:} In the regime of large $\tau$, Sinkhorn's linear dependence on $\tau$ much hinders its practicality, motivating \citep{sejourne2022faster} to alleviate this issue in the specific case of 1-D UOT. However, no algorithm for \emph{general} UOT problem that could achieve logarithmic dependence on $\tau$  had existed before GEM-UOT. From the practical perspective, $\tau$ can be large in certain applications \citep{Schiebinger_Optimal_2019} and, for a specific case of UOT penalized by squared $\ell_2$ norm\footnote{We recall that our paper considers UOT penalized by KL divergence, which is the most prevalent UOT variant (see more details in Section \ref{sec_intro}).}, be of the order of thousands \citep{blondel2018smooth}. 
    \item \textit{Logarithmic dependence on ${\varepsilon}^{-1}$:} GEM-UOT is the first algorithm in the literature that achieves logarithmic dependence on ${\varepsilon}^{-1}$. As the class of scaling algorithms, especially subsuming Sinkhorn, is known to suffer from  numerical instability \citep{deplaen2023unbalanced}, GEM-UOT thus emerges as a practical solver for the high-accuracy regime. 
\end{itemize}

\subsection{Relaxed UOT Problem and Practical Perspectives}
\label{sec:gemruot}

In this Section, we discuss the relaxed UOT problem where only the UOT distance is within interest, and propose a practical algorithm, termed GEM-RUOT, for this setting.  

\textbf{Motivation of Relaxed UOT Problem:}  Certain generative model applications \citep{Nguyen2021DistributionalSA}  recently have adopted either OT or UOT as loss metrics for  training, in which only the computation of the distance is required, while the transportation plan itself is not within interest. To this end, we consider the \textit{relaxed} UOT (RUOT) problem, where the goal is to approximate the UOT distance oblivious of the transportation plan. 

\textbf{Practicality:}
Gradient methods have been shown to achieve competitive complexity in  the OT literature \citep{guminov2021accelerated, dvurechensky2018computational}.  The heuristic  proposed by \citep{Yang_Scalable_2019} and based on gradient descent achieved favorable performance for solving the  UOT problem. 
Through our complexity analysis of GEM-UOT, we want to establish the preliminary understanding of the class of  gradient-based optimization when applied to UOT problem. Nevertheless, the convoluted function decomposition \eqref{primal-dual-form} of GEM-UOT may make it hard to implement in practice. The strong convexity constant $\mu$ (Lemma \ref{sparse_UOT_lipschitz_lemma}) can be small on certain real datasets, thereby hindering the empirical performance. We thus develop an easy-to-implement algorithm, called GEM-RUOT, for practical purposes that avoids the dependency on $\mu$ in complexity and is specialized for the RUOT problem. One interesting direction is to incorporate Adaptive Gradient Descent \citep{pmlr-v119-malitsky20a} that automatically adapts to  local geometry and is left for future work.


\textbf{Alternative Dual Formulation:} 
GEM-RUOT does not require intricate function decomposition and optimizes directly over the following  dual formulation \eqref{alt_dual_prob_adaptive1} in Lemma \ref{duality_theorem2}. 
Optimizing \eqref{alt_dual_prob_adaptive1} is equivalent to: 
\begin{equation}
     \min_{\uU, \vV \in \Br^{n}} h_a(\xX = (\uU, \vV)) := \frac{1}{4 \eta }\sum_{i,j=1}^{n} \max\{0, u_i+ v_j - C_{ij}\}^2 + \tau \left\langle  e^{-\uU/\tau} , \aA\right\rangle  + \tau \left\langle  e^{-\vV/\tau} , \bB\right\rangle, \label{dual_relaxed}
\end{equation}
which is locally smooth by the following Lemma \ref{rgem_convex_lipschitz_lemma}, whose proof is given in Appendix \ref{appen_lma16}. 


\begin{lemma}
\label{rgem_convex_lipschitz_lemma}
Let  $V_a = \{ \xX  = (\uU, \vV) \in \Br^{2n} :  \forall i,j \in [n],  \tau \log(\frac{2a_i}{\alpha+\beta}) \leq u_i \leq D,  \tau \log(\frac{2b_j}{\alpha+\beta}) \leq  v_j \leq D \}$. Then $h_a(\xX)$ is $L_a$-smooth and convex in $V_a$ with: 
$L_{a} = {\frac{\alpha+\beta}{\tau} + \frac{2\sqrt{n}}{\eta}}.$
\end{lemma}

GEM-RUOT adopts the  $\ell_2$ distance for its prox-function $\bar{P}(\xX_0, \xX) = \frac{1}{2}\|\xX - \xX_0\|_2^2$, and consequently uses the standard $\ell_2$ projection operator for the proximal mapping associated with the closed convex set $\Yy$ as 
$\bar{\mathcal{M}}_\Yy(\gG, \xX_0, \theta) := \argminE_{\xX = (\uU,\vV) \in \Yy} \{\langle \gG,\xX \rangle + \theta \bar{P}(\xX_0,\xX) \}.$
Using the convex version of GEM, GEM-RUOT optimizes directly over \eqref{dual_relaxed} and projects its iterates onto the closed convex set $V_a$ through the proximal mapping $\bar{\mathcal{M}}_{V_a}()$ to enforce the smoothness on the convex hull of all the iterates. Finally, it returns $F_a(\underline{\xX}^k)$ as an $\varepsilon$-approximation of the distance $\UOT_\KL(\aA, \bB)$ and a heuristic transportation plan $X^k$. The steps are summarized in Algorithm \ref{alg: acc-proximal-gd-uot-convex} and GEM-RUOT's complexity is  given in Theorem \ref{thm: acce_proximal_guarantee2}. 
\begin{theorem}
\label{thm: acce_proximal_guarantee2}
If $k \geq  \sqrt{\frac{12 L_a n D^2 }{\varepsilon}}$, then $|F_a(\underline{\xX}^k) -\UOT_\KL(\aA, \bB)|  \leq \varepsilon$.
Under the assumptions (A1-A3), the complexity of Algorithm~\ref{alg: acc-proximal-gd-uot-convex} is 
$O\left(  (\alpha + \beta) \cdot \frac{\tau n^{2.75}}{\varepsilon} \cdot \log(n(\alpha+\beta))\right)$.
\end{theorem}
\begin{proof}
Let $\xX^* = (\uU^*, \vV^*)$ be an optimal solution to \eqref{dual_relaxed} and thus \eqref{alt_dual_prob_adaptive1}.
Since $\|\xX^*\|_\infty \leq D$ by Lemma \ref{opt_sol_bounds} and Remark \ref{equi_sol_remark}, we have
$\bar{P}(\xX^0,\xX^*) = \frac{1}{2} \|\xX^*\|_2^2 \leq \frac{1}{2} n D^2$. 
From corollary 3.6 in \citep{Lan2018RandomGE}, we obtain:
\begin{align}
\nonumber
    0\leq h_a(\underline{\xX}^k) - h_a(\xX^*) = F_a(\xX^*) - F_a(\underline{\xX}^k) = &\UOT^{\eta}_\KL(\aA, \bB) - F_a(\underline{\xX}^k) \leq \frac{12L_a}{k(k+1)}\bar{P}(\xX^0,\xX^*)  \\
    \therefore 0\leq&\UOT^{\eta}_\KL(\aA, \bB) - F_a(\underline{\xX}^k) \leq  \frac{6L_a n D^2}{k(k+1)} \label{bobaboi1}
\end{align}
Using $\|X_f\|_1 \leq \frac{\alpha+\beta}{2}$ (Lemma \ref{opt_sol_bounds}) and $\eta = \frac{\varepsilon}{2R}=  \frac{2\varepsilon}{(\alpha+\beta)^2}$ set by the Algorithm \ref{alg: acc-proximal-gd-uot-convex}, we have:
\begin{align}
\nonumber
    \UOT^{\eta}_\KL(\aA, \bB) -\UOT_\KL(\aA, \bB) &=  g_\eta(X^\eta) - f(X_f)\\
    \nonumber
    &=  (f(X^\eta) - f(X_f))+\eta \|X^\eta\|_2^2  \geq 0\\
\nonumber
    \UOT^{\eta}_\KL(\aA, \bB) -\UOT_\KL(\aA, \bB) &= g_\eta(X^\eta) - g_\eta(X_f) + \eta \|X_f\|_2^2 \\
    \nonumber
    &\leq \eta \|X_f\|_2^2  \leq \eta \frac{(\alpha+\beta)^2}{4}  = \frac{\varepsilon}{2}\\
    \therefore 0 \leq  \UOT^{\eta}_\KL(\aA, \bB) -\UOT_\KL(\aA, \bB) &\leq  \frac{\varepsilon}{2} \label{bobaboi2}
\end{align}
Combining \eqref{bobaboi1}, \eqref{bobaboi2}, we obtain:
\begin{align}
\label{bobaboi3}
    |F_a(\underline{\xX}^k) - \UOT_\KL(\aA, \bB) | \leq \frac{\varepsilon}{2} + \frac{6L_a n D^2}{k(k+1)}
\end{align}
The RHS of \eqref{bobaboi3} is less than $\varepsilon$ if $k \geq \sqrt{\frac{12 L_a n D^2 }{\varepsilon}} = K_a$. 
Under the assumptions (A1-A3) and and $\eta = \frac{\varepsilon}{2R}=  \frac{2\varepsilon}{(\alpha+\beta)^2}$, we have:
\begin{align*}
    D &= O( \tau \log(n(\alpha+\beta))), \quad L_a = O\bigg(\frac{(\alpha+\beta)^2 \sqrt{n}}{\varepsilon}\bigg) \\
    \therefore K_a &= O\left(  (\alpha + \beta) \cdot \frac{\tau n^{0.75}}{\varepsilon} \cdot \log(n(\alpha+\beta))\right) 
\end{align*}
Since  it costs $O(n^2)$ per iteration to compute the gradient $h_a(\xX)$, the total complexity of Algorithm \ref{alg: acc-proximal-gd-uot-convex} is 
    $O\left(  (\alpha + \beta) \cdot \frac{\tau n^{2.75}}{\varepsilon} \cdot \log(n(\alpha+\beta))\right)$.    
\end{proof}

\section{Approximability of Standard OT}
\label{sec:approximation}


For balanced masses and $\tau \to \infty$, UOT problem reduces to standard OT. In this Section, we establish the very first characterization of the vanishing approximation error between UOT and OT in terms of $\tau$. To facilitate our discussion, we  formally define the standard OT problem in Section \ref{sec:standardOT}. 
We prove several some approximability properties of UOT in Section \ref{sec_help_approxi}, which are helpful for later deriving the main results. 
Then in Section \ref{approx_ot_solution}, we show that UOT's transportation plan converges to the marginal constraints of standard OT in the $\ell_1$ sense with the rate $O(\frac{n}{\tau})$. Beside transportation plan, we upperbound the difference in transport plan between  OT and  UOT distances by $O(\frac{n^2}{\tau})$ in Section \ref{sec_distance_error}. 




\subsection{Standard OT}
\label{sec:standardOT}

When the  masses lie in the probability simplex, i.e. $\aA, \bB \in \Delta^n:= \{\xX \in \Br^n: \|\xX\|_1 = 1 \}$,  the standard OT problem \citep{Kantorovich-1942-Translocation} can be formulated as: 
\begin{align}
\label{OT_func}
    \OT(\aA, \bB) = \min_{X \in \Pi(\aA, \bB)}  \left\langle C, X\right\rangle 
\end{align}
where $\Pi(\aA, \bB) = \{X \in \Br^{n\times n}: X \one_n = \aA, X^{\top} \one_n = \bB \}$.

\begin{definition}
For  $\varepsilon > 0$,  $X$ is an $\varepsilon$-approximation transportation plan of $\OT(\aA, \bB)$ if:
\begin{align}
    \left\langle C, X\right\rangle  \leq \left\langle C, X^{OT}\right\rangle +\varepsilon,
\end{align}
where $X^{OT} = \argminE_{X \in \Pi(\aA, \bB)} \left\langle C, X\right\rangle $ is the optimal transportation plan of the OT problem.
\end{definition}

\subsection{Approximability Properties of UOT}
\label{sec_help_approxi}

For any $\eta > 0$, consider the optimal solution $\xX^* = (\uU^*, \vV^*)$ to \eqref{alt_dual_prob_adaptive1}. 
Note that since $\aA, \bB \in \Delta^n$, we have $\alpha = \beta =1$. By Lemma \ref{opt_sol_bounds}, we obtain:
\begin{align}
\label{superbound1} 
\|X^\eta\|_1 \leq 1.
\end{align}
The following Lemmas are useful for the proofs of the main Theorem \ref{OT_solution_approx} and \ref{gem_ot1} later.

\begin{lemma}
\label{margin_bound_misc1}
The followings hold $\forall i, j \in [n]$:
\begin{align}
    u^*_i &\leq \|C\|_\infty + 2\eta \label{u_upper1}\\
    v^*_j &\leq \|C\|_\infty + 2\eta \label{v_upper1}
\end{align}
\end{lemma}
\begin{proof}
Firstly, we show that there must exist some $l\in[n]$ such that $v_l \geq 0$. Assume for the sake of contradiction that $ \forall j\in[n]: v_j < 0$. Then from \eqref{gradient_cond2}, we obtain $ \forall j\in[n]$: 
\begin{align*}
    \log(\sum_{k=1}^n X^\eta_{kj}) &> \log(b_j) \\(X^{\eta\top} \one_n)_j=\sum_{k=1}^n X^\eta_{kj} &> b_j\\
    \therefore \| X^\eta\|_1 = \sum_{j=1}^n  (X^{\eta\top} \one_n)_j &> \sum_{j=1}^n b_j = 1,
\end{align*}
where the last line contradicts \eqref{superbound1}.
From \eqref{superbound1}, \eqref{primal_sol} and taking any $v_l \geq 0$, we have $\forall i\in[n]$:
\begin{align*}
    &1 \geq X^\eta_{il} = \frac{1}{2\eta}\max\{0, u^*_i+ v^*_l - C_{il}\}\geq \frac{1}{2\eta} (u^*_i+ 0 - \|C\|_\infty) \implies u^*_i \leq \|C\|_\infty + 2\eta 
\end{align*}
Similarly, we can prove that $\forall j\in[n]$: $v^*_j \leq \|C\|_\infty + 2\eta $.

\end{proof}

\begin{lemma}
\label{margin_bound_misc2}
Define $\gamma= \|C\|_\infty + 2\eta$. Then the followings hold $\forall i, j \in [n]$:
\begin{align}
    a_i e^{-\gamma/\tau} &\leq \sum_{j=1}^n X^\eta_{ij} = ( X^\eta \one_n)_i \leq 1 - e^{-\gamma/\tau} (1-a_i) \label{margin_inter1} \\
     b_j e^{-\gamma/\tau} &\leq \sum_{i=1}^n X^\eta_{ij} = (X^{\eta\top} \one_n)_j \leq 1 - e^{-\gamma/\tau} (1-b_j) \label{margin_inter2}
\end{align}
\end{lemma}
\begin{proof}
From \eqref{gradient_cond1}, we have $\forall i\in[n]$: 
\begin{align}
\label{interbound1}
    ( X^\eta \one_n)_i = a_i e^{-u^*_i/\tau} \overset{\eqref{u_upper1}}{\geq}  a_i e^{-\gamma/\tau}
\end{align}
For the upper bound, we have $\forall i\in[n]$: 
\begin{align*}
    ( X^\eta \one_n)_i &= \| X^\eta\|_1 - \sum_{k \neq i} ( X^\eta \one_n)_k \overset{\eqref{superbound1}+\eqref{interbound1}}{\leq} 1 - \sum_{k \neq i} a_k e^{-\gamma/\tau} = 1 - e^{-\gamma/\tau} (1-a_i),
\end{align*}
which has  proved \eqref{margin_inter1}. Similarly, \eqref{margin_inter2} is obtained. 

\end{proof}

\begin{lemma}
\label{margin_bound1}
Define $\gamma= \|C\|_\infty + 2\eta$. Then the followings hold:
\begin{align}
    \| X^\eta\one_n - \aA \|_1 &\leq \frac{n\gamma}{\tau}  \label{margin_bound1.1}\\
    \|X^{\eta\top}\one_n - \bB \|_1 &\leq \frac{n\gamma}{\tau} \label{margin_bound1.2}
\end{align}
\end{lemma}
\begin{proof}
From \eqref{margin_inter1}, we have:
\begin{align*}
     -a_i(1- e^{-\gamma/\tau} )\leq ( X^\eta\one_n)_i - a_i &\leq (1-a_i)  (1-e^{-\gamma/\tau})\\
   \therefore   |( X^\eta\one_n)_i - a_i| &\leq \max\{a_i, 1-a_i\} (1-e^{-\gamma/\tau})\leq \gamma/\tau,
\end{align*}
where the last line follows from $\max\{a_i, 1-a_i\}\leq 1$ and  $e^x \geq x+1, \forall x\in \Br$. 
We  obtain that:
\begin{align*}
    \| X^\eta\one_n - \aA \|_1  = \sum_{i=1}^n |( X^\eta\one_n)_i - a_i| \leq \frac{n\gamma}{\tau},
\end{align*}
which has proved \eqref{margin_bound1.1}. Similarly, \eqref{margin_bound1.2} is obtained.
\end{proof}


\subsection{Approximation of OT Transportation Plan}
\label{approx_ot_solution}

We prove in the next Theorem the rate at which UOT's transportation plan converges to the marginal constraints of standard OT. 
\begin{theorem}
\label{OT_solution_approx}
For $\aA, \bB \in \Delta^n$,  the problem $\UOT_{\KL}(\aA, \bB)$ admits some optimal transportation plan $X_f = \argminE_{X \in \Br_{+}^{n \times n}} f(X)$   such that:
\begin{align}
\label{UOT_marginal_cons1}
     \| X_f\one_n - \aA \|_1 + \|X_f^{\top}\one_n - \bB \|_1 &\leq \frac{2n\|C\|_\infty}{\tau}.
\end{align}
\end{theorem}
\begin{proof}
    From Lemma \ref{opt_sol_bounds}, we know that $X^\eta$ is bounded by $\| X^\eta\|_1 \leq \frac{\alpha + \beta}{2}$. By the Bolzano–Weierstrass Theorem, there exists a convergent subsequence $\{ X^{\eta_t} \}_{t=1}^\infty$ where $\lim_{t \to \infty  } \eta_t = 0$. Define the limit of this subsequence as
    $X_f = \lim_{\eta_t \to 0} X^{\eta_t}$.
We first prove that $X_f$ is a transportation plan of the $\UOT_{\KL}(\aA, \bB)$ problem by showing that $X_f$ is a minimizer of $f$.
By the boundedness of $\| X^{\eta_t}\|_1 \leq \frac{\alpha + \beta}{2}$ (Lemma \ref{opt_sol_bounds}), we have $ \lim_{\eta_t \to 0} \eta_t \| X^{\eta_t} \|_2^2 = 0$ and obtain:
\begin{align}
     \lim_{\eta_t \to 0} g_{\eta_t}(X^{\eta_t}) = \lim_{\eta_t \to 0} [ f(X^{\eta_t}) + \eta_t \| X^{\eta_t} \|_2^2 ] = \lim_{\eta_t \to 0}  f(X^{\eta_t}) = f(\lim_{\eta_t \to 0} X^{\eta_t}) = f(X_f) \label{OT_solution_approx_proof1}
\end{align}
where the last equality is by the continuity of $f$. Now take any $\widetilde{X} \in \Br_+^{n \times n}$. we have:
\begin{align}
\nonumber
   g_{\eta_t}(X^{\eta_t}) &\leq  g_{\eta_t}(\widetilde{X}) = f(\widetilde{X}) + \eta_t \|\widetilde{X}\|_2^2\\
    \therefore \lim_{\eta_t \to 0} g_{\eta_t}(X^{\eta_t}) &\leq \lim_{\eta_t \to 0} g_{\eta_t}(\widetilde{X})  =  f(\widetilde{X}).  \label{OT_solution_approx_proof2}
\end{align}
Combining \eqref{OT_solution_approx_proof1} and \eqref{OT_solution_approx_proof2}, we have $ f(X_f) \leq f(\widetilde{X})$, thereby proving the minimality of $X_f$.
Thus, $X_f$ is a transportation plan of the $\UOT_{\KL}(\aA, \bB)$ problem. 
Now by Lemma \ref{margin_bound1},  we have: 
\begin{align*}
    \| X^{\eta_t}\one_n - \aA \|_1 + \|X^{\eta_t\top}\one_n - \bB \|_1 &\leq \frac{2n(\|C\|_\infty + 2\eta_t)}{\tau}
\end{align*}
Taking $\eta_t \to 0$ and noting that $X_f = \lim_{\eta_t \to 0} X^{\eta_t}$, we obtain the desired inequality \eqref{UOT_marginal_cons1}.
\end{proof}


The quantity $\| X\one_n - \aA \|_1 + \|X^{\top}\one_n - \bB \|_1 $ measures the closeness of $X$ to $\Pi(\aA, \bB)$ in $\ell_1$ distance and is also used as the stopping criteria in algorithms solving standard OT \citep{Altschuler-2018-Approximating, dvurechensky2018computational}. 
Therefore, Theorem \ref{OT_solution_approx} characterizes the rate $O(\frac{n}{\tau})$ by which the transportation plan of $\UOT_{\KL}(\aA. \bB)$ converges to the marginal constraints $\Pi(\aA, \bB)$ of $\OT(\aA. \bB)$.
Utilizing this observation, we present the algorithm  GEM-OT (Algorithm \ref{alg: gem-ot}), which solves the  $\UOT_{\KL}(\aA. \bB)$ problem via GEM-UOT with fine tuned $\tau$ and finally projects the UOT solution onto $\Pi(\aA, \bB)$ (via Algorithm \ref{alg: proj}), to find an $\varepsilon$-approximate solution to the standard OT problem $\OT(\aA, \bB)$. The following Theorem characterizes the complexity of GEM-OT.

\begin{minipage}[t]{0.46\textwidth}
\begin{algorithm}[H]  
\caption{GEM-OT}
\label{alg: gem-ot}
\begin{algorithmic}[1]
 \STATE {\bfseries Input:} $C, \aA, \bB, \varepsilon$
 \STATE $\varepsilon' = \varepsilon/16$
 \STATE $\eta  = \frac{2\varepsilon'}{(\alpha+\beta)^2} = \varepsilon'/{2}$
 \STATE $ \gamma = \|C\|_\infty +\eta$
 \STATE $\tau =   \frac{16\| C\|_\infty n \gamma}{\varepsilon}$
\STATE $\bar{X}=$ {GEM-UOT}$(C, \aA, \bB, \varepsilon', \tau, \eta)$ 
\STATE $Y =$ {PROJ}$(\bar{X}, \aA, \bB)$
\STATE {\bfseries Return} $Y$\\
\end{algorithmic}
\end{algorithm}
\end{minipage}
\hfill
\begin{minipage}[t]{0.46\textwidth}
\begin{algorithm}[H] 
\caption{PROJ  \cite[Algorithm 2]{altschuler2017near}}
\label{alg: proj}
\begin{algorithmic}[1]
\STATE {{\bf Input:} $X \in \Br^{n\times n}, \aA \in \Br^{n}, \bB\in \Br^{n}$}
\STATE $P \leftarrow diag(\xX)$ with $x_i = \frac{a_i}{(X\one_n)_i} \wedge 1$
\STATE $X' \leftarrow P X$
\STATE $Q \leftarrow diag(\yY)$ with $y_j = \frac{b_j}{(X'^{\top}\one_n)_j} \wedge 1$
\STATE $X'' \leftarrow X' Q$
 \STATE $err_r \leftarrow \aA - X'' \one_n $,  $err_c \leftarrow \bB - X''^{\top} \one_n $
 \STATE {\bfseries Return} $Y \leftarrow X'' + err_r err_c^T/\|err_r\|_1$\\
\end{algorithmic}
\end{algorithm}
\end{minipage}

\begin{theorem}
\label{gem_ot1}
Algorithm \ref{alg: gem-ot} outputs $Y\in \Pi(\aA, \bB)$ such that:  
$$\left\langle C, Y\right\rangle  \leq \left\langle C, X^{OT}\right\rangle +\varepsilon.$$
In other words, $Y$ is an $\varepsilon$-approximate transportation plan of $\OT(\aA, \bB)$. 
Under assumptions (A1-A2), the complexity of Algorithm \ref{alg: gem-ot} is 
$\widetilde{O}\big(  {\kappa}  {n^2} \big)$.
\end{theorem}
\begin{proof}
We obtain via simple algebra and Holder inequality that:
\begin{align}
\nonumber 
    \left\langle C, Y\right\rangle -  \left\langle C, X^{OT}\right\rangle  &= \left\langle C,  \bar{X} - X^\eta \right\rangle + \left\langle C, Y - \bar{X}\right\rangle +    \left\langle C,   X^\eta - X^{OT} \right\rangle \\
    &\leq  \| C\|_\infty\|  \bar{X} - X^\eta \|_1 +\| C\|_\infty\| Y - \bar{X}\|_1 +   \left\langle C,   X^\eta - X^{OT} \right\rangle  \label{breakdown1}
\end{align}
We proceed to upper-bound each of the three terms in the RHS of \eqref{breakdown1}.
We note that  $\bar{X}$ is the output of GEM-UOT fed with input error $\varepsilon/16$, and GEM-UOT (in the context of  Algorithm \ref{alg: gem-ot})  thus sets $\eta = \frac{2\varepsilon/16}{(\alpha+\beta)^2} = \frac{\varepsilon}{32}$.
The first term can be  bounded by Lemma \ref{lma: acce_proximal_guarantee1} as:
\begin{align}
    \| C\|_\infty\|  \bar{X} - X^\eta \|_1 \leq  \| C\|_\infty \frac{\varepsilon/16}{2L_1} \leq \frac{\varepsilon}{8} \label{term1}
\end{align}
where $L_1 \geq \|C\|_\infty$ holds directly from the definition of $L_1$. 
For the second term, we obtain from Lemma \ref{proj_algo} that:
\begin{align}
\nonumber 
    \| C\|_\infty\| Y - \bar{X}\|_1 &\leq  2 \| C\|_\infty [\|\bar{X} \one_n - \aA \|_1 + \|\bar{X}^{\top} \one_n - \bB \|_1]\\ 
    \nonumber
    &\leq 2 \| C\|_\infty [\| X^\eta \one_n - \aA \|_1 + \|{X}^{\eta\top} \one_n - \bB \|_1] \\
    \nonumber 
    &+ 2 \| C\|_\infty [\|( X^\eta - \bar{X}) \one_n  \|_1 + \|({X}^{\eta\top} - \bar{X}^{\top} )\one_n  \|_1]\\ 
    &\leq \frac{4\| C\|_\infty n \gamma}{\tau} +4 \| C\|_\infty\|\bar{X} - X^\eta\|_1   \label{explain1} \\
    & \leq \frac{\varepsilon}{4} +  \frac{4\varepsilon}{8} = \frac{3\varepsilon}{4} \label{term2}
\end{align}
where \eqref{explain1} follows Lemma \ref{margin_bound1} with $\gamma = \|C\|_\infty + 2\eta$ and triangular inequality, and \eqref{term2} is by the definition of $\tau$ and \eqref{term1}. 
Since $X^{OT}\in \Pi(\aA, \bB)$, we have: $\KL(X^{OT} \one_n||\aA) = \KL(X^{OT \top} \one_n||\bB) =0 $ and  $\eta \|X^{OT}\|_2^2 \leq \eta \|X^{OT}\|_1^2 = \eta$. We thus obtain that:
\begin{align}
\label{minibound1}
    g_\eta(X^{OT}) = \left\langle C,  X^{OT} \right\rangle + \eta \| X^{OT}\|_2^2  \leq \left\langle C,  X^{OT} \right\rangle + \eta
\end{align}
Since the $\KL$ divergence and $\ell_2$-norm are non-negative, we have:
\begin{align}
\label{minibound2}
    g_\eta( X^\eta) \geq  \left\langle C,  X^{\eta} \right\rangle. 
\end{align}
We recall that $ X^\eta  = \argminE_{X \in \Br_{+}^{n \times n}} g_\eta(X)  $, which implies:
\begin{align}
\label{minibound3}
    g_\eta(X^{OT}) \geq g_\eta( X^\eta).
\end{align}
Combining \eqref{minibound1}, \eqref{minibound2} and \eqref{minibound3}, we upper-bound the third term in the RHS of \eqref{breakdown1}:
\begin{align}
    \left\langle C,  X^{\eta}- X^{OT} \right\rangle \leq \eta = \frac{\varepsilon}{32} \label{term3}.
\end{align}
Plugging \eqref{term1}, \eqref{term2} and \eqref{term3} into \eqref{breakdown1}, we have:
\begin{align*}
     \left\langle C, Y\right\rangle -  \left\langle C, X^{OT}\right\rangle \leq \frac{\varepsilon}{8} + \frac{3\varepsilon}{4} + \frac{\varepsilon}{32} < \varepsilon. 
\end{align*}

The complexity of GEM-OT is the total of $O(n^2)$ for the PROJ() algorithm and the complexity of GEM-UOT for the specified $\tau$ in the algorithm. By Corollary \ref{acceGD_complexity1}, which establishes the complexity of GEM-UOT under the assumptions (A1-A3), we conclude that the complexity of GEM-OT under (A1-A2) (where (A3) is naturally satisfied by our choice of $\tau$) is
$\widetilde{O}\big(  {\kappa}  {n^2}  \big)$. 
\end{proof}


\begin{remark}
The best-known complexities in the literature of OT are respectively $\widetilde{O}(n^{2.5})$ \citep{6979027} for finding exact solution, and   $\widetilde{O}(\frac{n^2}{\varepsilon})$ \citep{pmlr-v206-luo23a, yilingxie2023accelerated_pdasgd} for finding the $\varepsilon$-approximate solution. We highlight that the complexity  in  Theorem \ref{gem_ot1} is the first  to achieve logarithmic dependence on $\varepsilon^{-1}$ to obtain the  $\varepsilon$-approximate solution, while still enjoying better dependence on $n$ than the exact solver \citep{6979027}. As the prevalent entropy-regularized methods to approximate OT \citep{pmlr-v206-luo23a, altschuler2017near, dvurechensky2018computational} have been known to suffer from  numerical instability and blurring issues \citep{zhou2023efficient}, GEM-OT can emerge as a practical solver for OT  specialized in the high-accuracy regime \citep{dong2020study}.
\end{remark}
Theorem \ref{OT_solution_approx} would elucidate the tuning  of $\tau$, which has been a heuristic process in all prior work, such that the transportation plan of $\UOT_{\KL}(\aA, \bB)$ well respects the original structure of that of $\OT(\aA, \bB)$. 
Furthermore, we provide a novel recipe for  OT retrieval from UOT, which combines Theorem \ref{OT_solution_approx} with a low-cost post-process projection step and can be integrated with any UOT solver beside GEM-UOT as considered in this paper.


\subsection{Approximation Error between UOT and OT Distances} 
\label{sec_distance_error}


As discussed in Section \ref{sec:open_problems}, while a lot of work have used UOT with small $\tau = 1$ as a relaxed variant of OT \citep{Fatras2021UnbalancedMO, balaji2020robust, le2021robust}, we could empirically illustrate noticeably large approximation error between the UOT and OT distances on real dataset. 
Such deviation can be  detrimental to target applications including generative modelling  \citep{Fatras2021UnbalancedMO} which recently has adopted  UOT in place of OT as a loss metrics for  training. This further necessitates the understanding of approximation error between UOT and OT in the sense of purely transport distance.  
We next establish an upper bound on the approximation error between UOT and OT distances in  Theorem \ref{thrm: distance_error}. 



\begin{theorem}
\label{thrm: distance_error}
For $\aA, \bB \in \Delta^n$, we define $M = \log(2) \|C\|_\infty^2 \big(n+3\kappa \big)^2 + 2n \|C\|_\infty^2$ and have the following bound: 
\begin{align}
    0 \leq \OT(\aA, \bB) - \UOT_{\KL}(\aA, \bB) \leq \frac{M}{\tau}.
\end{align}
\end{theorem}
\begin{proof}
The lower bound is straightforward. Since $X^{OT}$ (the transportation plan of $\OT(\aA, \bB)$) is a feasible solution to the optimization problem $\UOT_\KL(\aA, \bB)$ and $\KL(X^{OT}\one_n||\aA) = \KL(X^{OT\top}\one_n||\bB) = 0$ (as $X^{OT} \in \Pi(\aA, \bB)$), we obtain that:
\begin{align*}
    \OT(\aA, \bB) =\left\langle C, X^{OT}\right\rangle = f(X^{OT}) \geq f(X_f) = \UOT_\KL(\aA, \bB)
\end{align*}
We proceed to prove the upper bound of the error approximation. 
We first define a variant of $\UOT_{\KL}(\aA, \bB)$ problem where the $\KL()$ divergence is replaced by squared $\ell_2$ norm:
\begin{equation}
\label{UOT_l2}
    \UOT_{\ell_2}(\aA, \bB) =  \min_{X \in \Br_{+}^{n \times n}} \big\{ f_{l_2}(X) := \left\langle C, X\right\rangle   + \frac{\tau}{2\log(2)} \|X \one_{n} - \aA\|_2^2 + \frac{\tau}{2\log(2)} \|X^{\top} \one_{n} - \bB\|_2^2  \big\},
\end{equation}
By \citep[Theorem 2]{blondel2018smooth}, for $M_1 = \log(2) \|C\|_\infty^2 \big(n+3\kappa \big)^2$, we have:
\begin{align}
\label{l2_error_theorem}
     0 \leq \OT(\aA, \bB) -\UOT_{\ell_2}(\aA, \bB) \leq \frac{M_1}{\tau}.
\end{align}

\begin{remark}
We  highlight that deriving the  approximation error for $\UOT_{\KL}(\aA, \bB)$ is harder than $\UOT_{\ell_2}(\aA, \bB)$. First, KL divergence is non-smooth, which is  stated in ~\citep{blondel2018smooth} as the main reason for their choice of squared $\ell_2$ norm instead of KL to relax the marginal constraints. Second, 
the proof in \citep{blondel2018smooth} is based on interpreting the dual of $\UOT_{\ell_2}(\aA, \bB)$  as the sum of the dual of $\OT(\aA, \bB)$  and the a squared $\ell_2$-norm regularization on the dual variables, so the approximation error is derived by  simply bounding the dual variables. 
 \end{remark}
Recall that:
\begin{align*}
    \UOT^{\eta}_{\KL}(\aA, &\bB) =  \min_{X \in \Br_{+}^{n \times n}} \{ g_\eta(X) := \left\langle C, X\right\rangle + \eta \| X\|_2^2+ \tau \KL(X \one_{n} || \aA) + \tau \KL(X^{\top} \one_{n} || \bB)\},
\end{align*}
and $ X^\eta  = \argminE_{X \in \Br_{+}^{n \times n}} g_\eta(X)  $ is the solution to $ \UOT^{\eta}_{\KL}(\aA, \bB)$.
We further define the problem $ \overline{\UOT}^{\eta}_{\KL}(\aA, \bB)$ as follows:
\begin{equation*}
    \overline{\UOT}^{\eta}_{\KL}(\aA, \bB) =  \min_{X \in \Br_{+}^{n \times n}, \|X\|_1=1} \{ g_\eta(X) := \left\langle C, X\right\rangle + \eta \| X\|_2^2+ \tau \KL(X \one_{n} || \aA) + \tau \KL(X^{\top} \one_{n} || \bB)\},
\end{equation*}
which optimizes over the same objective cost as $\UOT^{\eta}_{\KL}(\aA, \bB)$ yet with the additional constraint that $X$ lies in the probability simplex, i.e. $\|X\|_1=1$. We let   the solution to $\overline{\UOT}^{\eta}_{\KL}(\aA, \bB)$ be $Z^\eta = \argminE_{X \in \Br_{+}^{n \times n}, \|X\|_1=1}  g_\eta(X)$. Since $\|Z^\eta\|_1 =1 $, we have $Z^\eta \one_n, Z^{\eta \top} \one_n \in \Delta^n$. Then Pinsker inequality gives:
\begin{align*}
    \KL(Z^\eta \one_{n} || \aA) &\geq \frac{1}{2\log(2)} \|Z^\eta \one_{n} - \aA\|_1^2 \geq  \frac{1}{2\log(2)} \|Z^\eta \one_{n} - \aA\|_2^2 \\
    \KL(Z^{\eta\top} \one_{n} || \bB) &\geq \frac{1}{2\log(2)} \|Z^{\eta\top} \one_{n} - \bB\|_1^2 \geq  \frac{1}{2\log(2)} \|Z^{\eta\top} \one_{n} - \bB\|_2^2 
\end{align*}
Consequently, we obtain that:
\begin{align}
    \overline{\UOT}^{\eta}_{\KL}(\aA, \bB)= g_\eta(Z^\eta) \geq f_{\ell_2}(Z^\eta) \geq  \UOT_{\ell_2}(\aA, \bB) \label{KL_to_l2}
\end{align}
For any $\eta > 0$, we consider $Y^\eta = PROJ(X^\eta, \aA, \bB)$ as the projection of $X^\eta$ onto $\Pi(\aA, \bB)$ via Algorithm \ref{alg: proj}. Then $\|Y^\eta\|_1= \aA^{\top} \one_n = 1$. 
By Lemma \ref{proj_algo} and \ref{margin_bound1}, we have: 
\begin{align}
    \|Y^\eta - X^\eta\|_1 &\leq 2 [\|X^\eta \one_n - \aA \|_1 + \|X^{\eta\top} \one_n - \bB \|_1] \leq \frac{2n(\|C\|_\infty + 2\eta)}{\tau}.  \label{Y_X_bound}
\end{align}
Since $Y^\eta$ is a feasible solution to  the optimization problem  $\overline{\UOT}^{\eta}_{\KL}(\aA, \bB)$, we obtain that:
\begin{align}
    g_\eta(Y^\eta) &\geq \overline{\UOT}^{\eta}_{\KL}(\aA, \bB), \label{Y_vs_UOT_simplex} 
\end{align}
Note that $\KL(Y^{\eta}\one_n||\aA) = \KL(Y^{\eta\top}\one_n||\bB) = 0$ as $Y^{\eta} \in \Pi(\aA, \bB)$. We thus can write
    $g_\eta(Y^\eta) =  \left\langle C, Y^\eta\right\rangle + \eta \|Y^\eta\|_2^2$
and obtain via Holder inequality that:
\begin{align*}
    g_\eta(Y^\eta) - \UOT^{\eta}_{\KL}(\aA, \bB)   &=  g_\eta(Y^\eta) - g_\eta(X^\eta) \\ 
    &=  \left\langle C, Y^\eta - X^\eta \right\rangle + \eta (\|Y^\eta\|_2^2 -  \|X^\eta\|_2^2) - \tau \KL(X^\eta \one_{n} || \aA) - \tau \KL(X^{\eta\top} \one_{n} || \bB)\\
    &\leq \|C\|_\infty \| Y^\eta - X^\eta \|_1 + \eta  \|Y^\eta\|_2^2 \overset{\eqref{Y_X_bound}}{\leq} \frac{2n \|C\|_\infty (\|C\|_\infty + 2\eta)}{\tau} + \eta
\end{align*}
where for the last inequality we note that $ \|Y^\eta\|_2 \leq \|Y^\eta\|_1 = 1$. The above is equivalent to:
\begin{align}
     \UOT^{\eta}_{\KL}(\aA, \bB) &\geq  g_\eta(Y^\eta) - \frac{2n \|C\|_\infty (\|C\|_\infty + 2\eta)}{\tau} - \eta.\label{UOT_to_Y}
\end{align}
Now combining \eqref{KL_to_l2}, \eqref{Y_vs_UOT_simplex}  and \eqref{UOT_to_Y}, we have $\forall \eta >0$:
\begin{align*}
    \UOT^{\eta}_{\KL}(\aA, \bB) \geq \UOT_{\ell_2}(\aA, \bB) - \frac{2n \|C\|_\infty (\|C\|_\infty + 2\eta)}{\tau} - \eta
\end{align*}
Taking the limit $\eta \to 0$ on both sides and noting that $\lim_{\eta \to 0} \UOT^{\eta}_{\KL}(\aA, \bB) = \UOT_{\KL}(\aA, \bB)$, we obtain:
\begin{align*}
    \UOT_{\KL}(\aA, \bB) &\geq \UOT_{\ell_2}(\aA, \bB) - \frac{2n \|C\|_\infty^2 }{\tau} \\
    \OT(\aA, \bB) - \UOT_{\KL}(\aA, \bB) &\leq  \OT(\aA, \bB)- \UOT_{\ell_2}(\aA, \bB) + \frac{2n \|C\|_\infty^2 }{\tau} 
\end{align*}
Combining with \eqref{l2_error_theorem}, we get the desired bound 
$\OT(\aA, \bB) - \UOT_{\KL}(\aA, \bB) \leq \frac{M}{\tau}$
where $M = M_1 + 2n \|C\|_\infty^2$. 
\end{proof}

\section{Experiments}
\label{sec_experiment}

In this Section, we  empirically verify GEM-UOT's performance and sparsity, and our theories on approximation error. To compute the ground truth value, we use the convex programming package \textbf{cvxpy} to find the exact UOT plan \citep{agrawal2019rewriting}. Unless specified otherwise, we always set $\eta = \frac{\varepsilon}{2R} = \frac{2\varepsilon}{(\alpha+\beta)^2}$ (according to Theorem \ref{thm: acce_proximal_guarantee1}).



\subsection{Synthetic Data}

%

We set $n = 200, \tau = 55,  \alpha = 4, \beta = 5$. Then $a_i$'s are drawn from uniform distribution and rescaled to ensure $\alpha = 4$, while $b_i$'s are drawn from normal distribution with $\sigma = 0.1$. Entries of $C$ are drawn uniformly from $[10^{-1}, 1]$. For both GEM-UOT and Sinkhorn, we vary $\varepsilon= 1 \to 10^{-4}$ to evaluate their time complexities in Figure 
\ref{fig_001} 
and test the $\tau$ dependency in Figure \ref{fig_002}. Both experiments verify GEM-UOT's better theoretical dependence on $\tau$ and $\varepsilon$.\\
\begin{minipage}[t]{0.46\textwidth}
\begin{figure}[H]  
\centering
\includegraphics[height=0.14\textheight,width=1.07\textwidth]{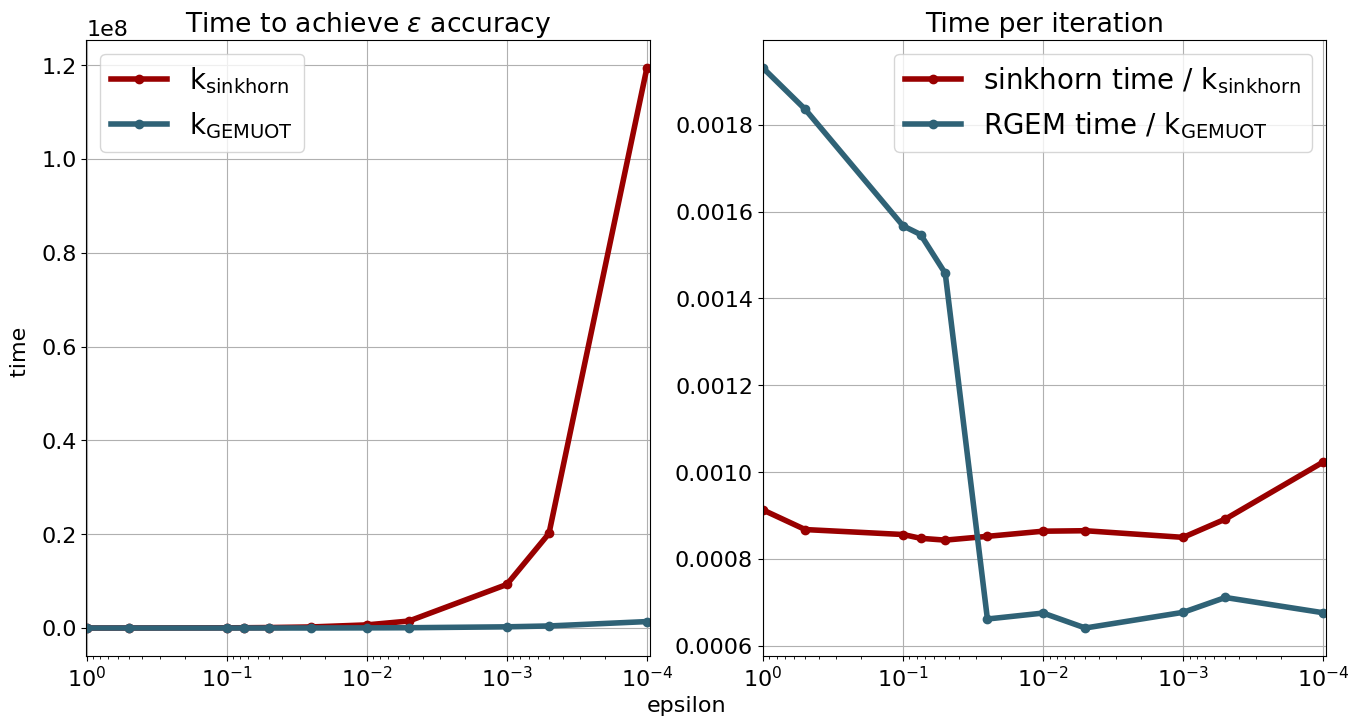}
\caption{Comparison in  number of iterations (left) and per-iteration cost (right) between GEM-UOT  and Sinkhorn  on  synthetic data for $\epsilon=1 \to 10^{-4}$.}
\label{fig_001}
\end{figure}
\end{minipage}
\hfill
\begin{minipage}[t]{0.46\textwidth}
\begin{figure}[H]   
\centering
\includegraphics[height=0.14\textheight,width=1.07\textwidth]{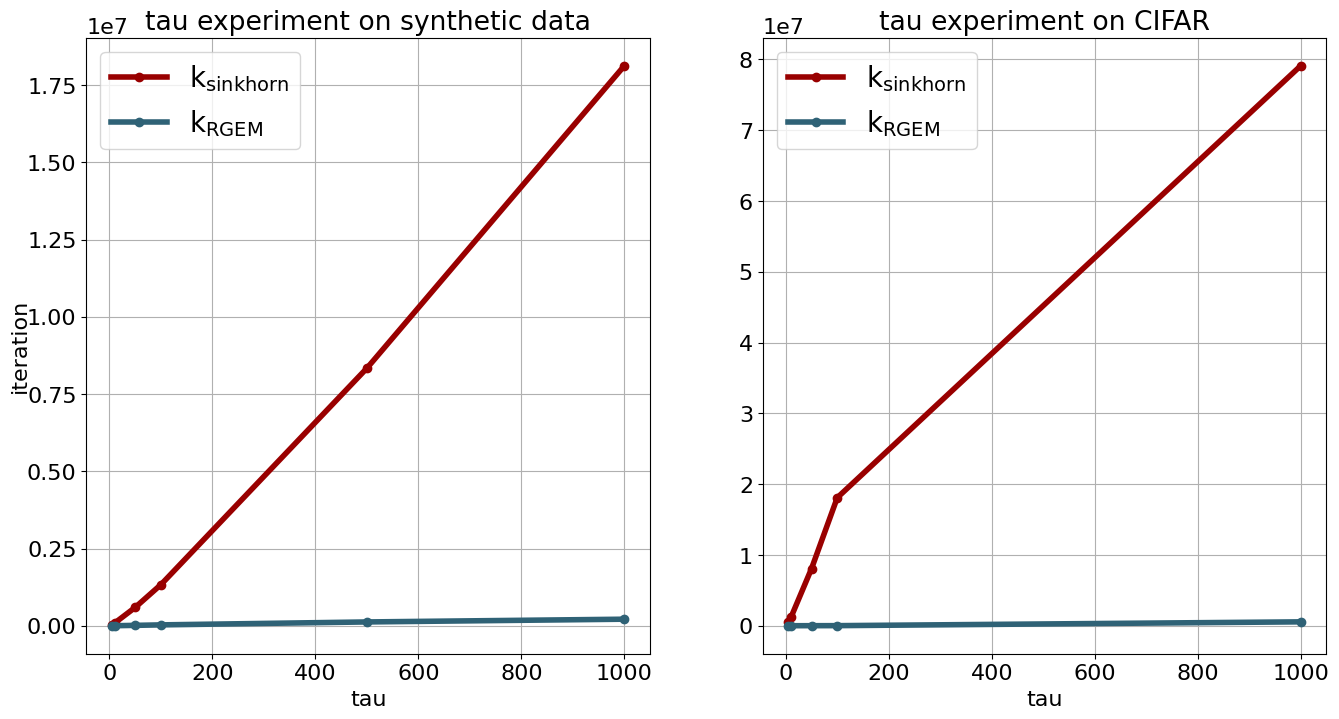}
\caption{Scalability of $\tau$ (for  $\varepsilon= 10^{-2}$) on synthetic data  and CIFAR-10. Sinkhorn  scales linearly while GEM-UOT scales logarithmically in $\tau$.}
\label{fig_002}
\end{figure}
\end{minipage}
\\

\begin{figure}[H] 
\centering
\includegraphics[height=0.193\textheight ,width=0.4\textwidth]{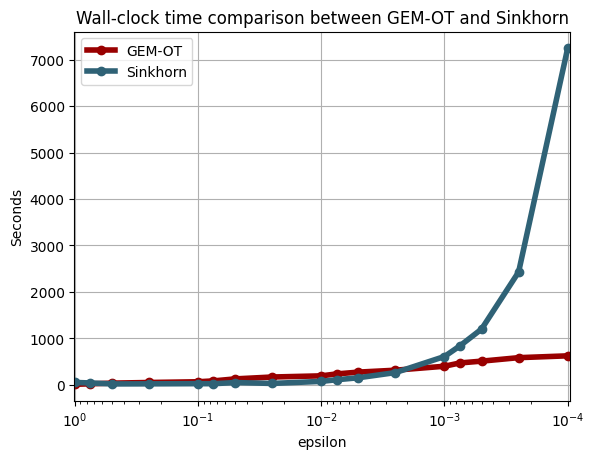}
\caption{Wall-clock time comparison between GEM-OT and Sinkhorn on synthetic dataset.}
\label{fig_gem_ot}
\end{figure}

Next, to empirically evaluate the efficiency of GEM-OT for approximating standard OT problem, we consider the  setting of balanced masses and large $\tau$ with $n =200, \alpha = \beta = 1$ and $\tau = 500$. We compare GEM-OT with Sinkhorn \citep{cuturi2013sinkhorn} in wall-clock time for varying $\varepsilon= 1 \to 10^{-4}$  in Figure \ref{fig_gem_ot}. 


\subsection{Real Datasets: CIFAR-10 and Fashion-MNIST}
\label{sec:real_data_exp}
To practically validate our algorithms, we compare the GEM-UOT with Sinkhorn on CIFAR-10 dataset \citep{Krizhevsky09learningmultiple}. A pair of flattened images corresponds to the marginals $\aA, \bB$, whereby the cost matrix $C$ is the matrix of $\ell_1$ distances between pixel locations. This is also the setting  considered in \citep{pham2020unbalanced, dvurechensky2018computational}.
 We plot the results in Figure \ref{fig_003}, which demonstrates GEM-UOT' superior performance. Additional experiments on Fashion-MNIST dataset that further cover GEM-RUOT are presented in Appendix \ref{sec:mnist1}.
We highlight that GEM-UOT only has logarithmic dependence on $\tau$, improving over Sinkhorn's linear dependence on  $\tau$ which had been its major bottleneck in the regime of large $\tau$ \citep{sejourne2022faster}. 
Experiment on the scalability of $\tau$ in Figure \ref{fig_002}, which is averaged over 10 randomly
chosen image pairs, empirically reasserts our favorable dependence on $\tau$ on the CIFAR-10 dataset.
In addition, we  evaluate the behaviour of GEM-UOT for different levels of $\varepsilon$ and thus $\eta$ in Figure \ref{fig_eta}.\\
\begin{minipage}[t]{0.46\textwidth}
\begin{figure}[H]
\centering
\includegraphics[height=0.2\textheight,width=1.0\textwidth]{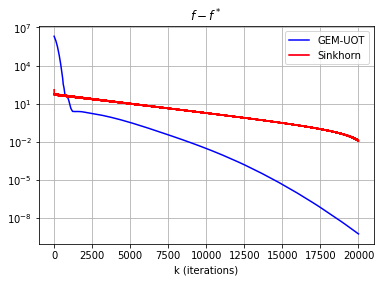}
\caption{Primal gap $f-f^*$ of GEM-UOT and Sinkhorn on  CIFAR-10.}
\label{fig_003}
\end{figure}
\end{minipage} 
\hfill
\begin{minipage}[t]{0.46\textwidth}
\begin{figure}[H] 
\centering
\includegraphics[height=0.193\textheight ,width=1.0\textwidth]{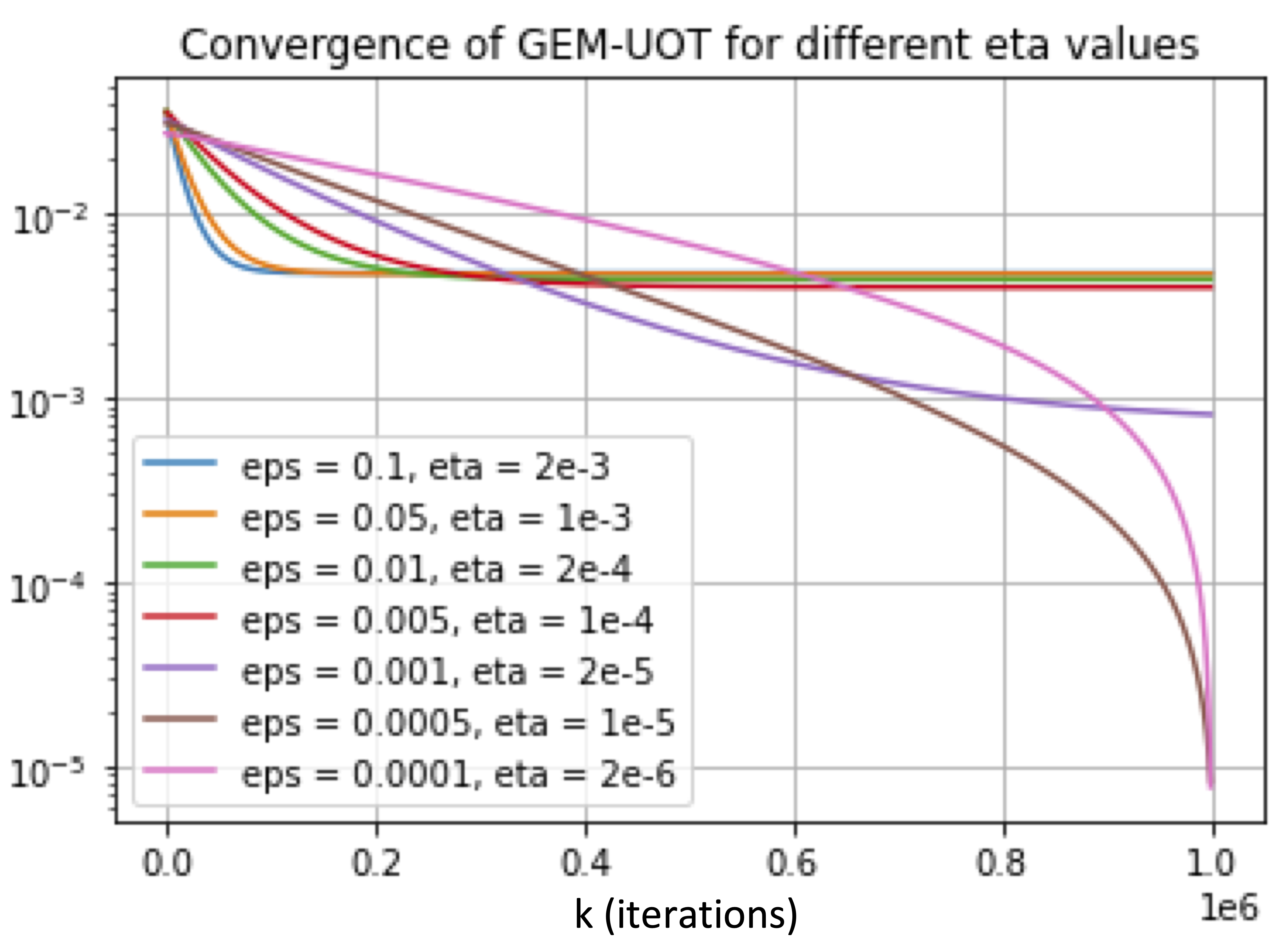}
\caption{Algorithmic dependence on $\eta$ and the corresponding accuracy achieved by the algorithm on  CIFAR-10  for $10^6$ iterations.}
\label{fig_eta}
\end{figure}
\end{minipage}



.

\subsection{Sparsity of GEM-UOT}
\label{sparse_sec}
On two real datasets Fashion-MNIST and CIFAR 10, we empirically illustrate the sparseness of the transportation plans produced by GEM-UOT in Figure \ref{fig_005} and Figure \ref{fig_007}, while Sinkhorn produces transportation plan with $0\%$ sparsity in those experiments. This result is consistent with \citep{blondel2018smooth}, where entropic regularization enforces strictly positive and dense transportation plan and squared-$\ell_2$ induces sparse transportation plan. To illustrate the effectiveness of $\ell_2$ as the sparse regularization, in Figure \ref{fig_compare1} we compare the heat maps of the optimal transportation plans of  squared-$\ell_2$ UOT  and original UOT  on CIFAR 10. We set the sparsity threshold of $10^{-2}$ and report the sparsity results in the following figures.\\
\begin{minipage}[b]{0.25\textwidth}
\begin{figure}[H]  
\centering
\includegraphics[width=1.1\textwidth]{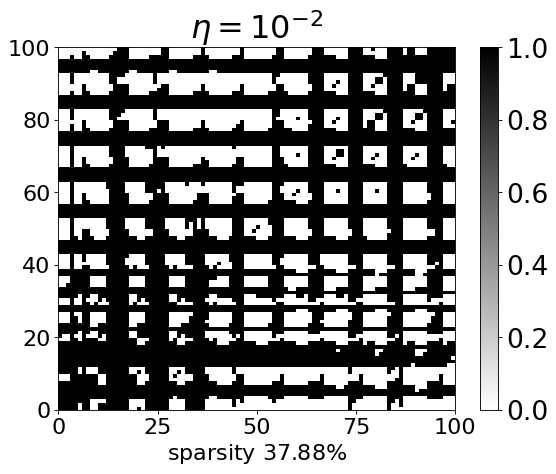} 
\caption{$37.88\%$ sparsity on Fashion-MNIST.}
\label{fig_005}
\end{figure}
\end{minipage}
\hfill
\begin{minipage}[b]{0.25\textwidth}
\begin{figure}[H]   
\centering
\includegraphics[width=1.1\textwidth]{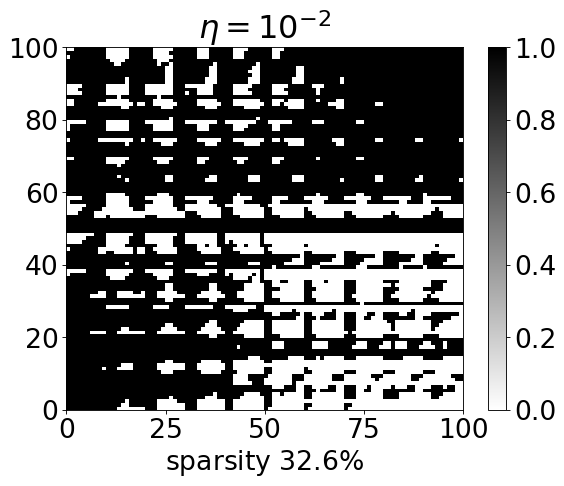}
\caption{$32.6\%$ sparsity on CIFAR-10.}
\label{fig_007}
\end{figure}
\end{minipage}
\hfill
\begin{minipage}[b]{0.46\textwidth}
\begin{figure}[H]   
\centering
\includegraphics[width=0.9\textwidth]{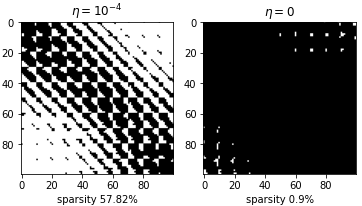}
\caption{ Sparsity of squared-$\ell_2$  UOT plan (left) versus original UOT plan (right)  on CIFAR 10.}
\label{fig_compare1}
\end{figure}
\end{minipage}



\subsection{Color Transfer}
Color transfer is a prominent application where sparse transportation plan is of special interest \citep{blondel2018smooth}. This experiment shows that  squared-$\ell_2$ UOT (i.e. sparse UOT) solved by GEM-UOT results in a sparse transportation plan and is thus congruent with this application, while entropic UOT solved by Sinkhorn results in strictly positive and dense transportation plan, hindering its adoption. 
In particular, we conduct the UOT  between the histograms of two images. Given the source image $X$ and target image $Y$ both of size $l \times h$, we present their pixels in RGB space. Quantizing the images down to $n$ colors, we obtain the color centroids for the source and target images respectively as $S_{source} = \{x_1, x_2, ..., x_n\}$ and $S_{target} = \{y_1, y_2,..., y_n\}$, where $x_i, y_i \in \Br^3$. Then the source image's color histogram $\aA \in \Delta^n$ is computed as the empirical distribution of pixels over the $n$ centroids, i.e. $a_k = \frac{1}{l  h}\cdot \big(\sum_{i=1}^l \sum_{j=1}^h \mathbbm{1}\{X_{ij} \in x_k\}\big)$. The target image's color histogram $\bB \in \Delta^n$ is computed similarly. 
We compute the transportation plan between the two images via sparse UOT solved by GEM-UOT and entropic UOT solved by Sinkhorn, and report their qualities of color transfer in Figure \ref{color_transfer_quality}, and the sparsity in Figure \ref{gem-uot-sparsity-color-transfer} and Figure \ref{sinkhorn-sparsity-color-transfer}. To demonstrate how entropic regularized UOT tend to generate dense mappings, we set the sparsity threshold of $10^{-2}$ for Figure \ref{sinkhorn-sparsity-color-transfer} while the sparsity threshold for Figure \ref{gem-uot-sparsity-color-transfer} is $0$. (Note that due to the entropic regularization, Sinkhorn has $0\%$ sparsity  under such  threshold.)
\\
\begin{minipage}[b]{0.53\textwidth}
\begin{figure}[H]
\centering
\includegraphics[width=1\textwidth]{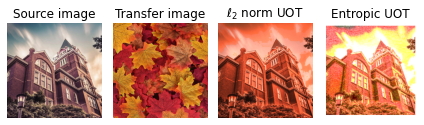}
\caption{Quality of color transfer by  GEM-UOT (solving squared-$\ell_2$  UOT) and Sinkhorn (solving entropic UOT).}
\label{color_transfer_quality}
\end{figure}
\end{minipage}
\hfill
\begin{minipage}[b]{0.21\textwidth}
\begin{figure}[H]
\includegraphics[width=1.02\textwidth]{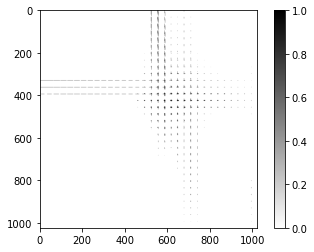}
\caption{GEM-UOT  gets $99.4\%$ sparsity.}
\label{gem-uot-sparsity-color-transfer}
\end{figure}
\end{minipage}
\hfill
\begin{minipage}[b]{0.204\textwidth}
\begin{figure}[H]
\includegraphics[width=1.03\textwidth]{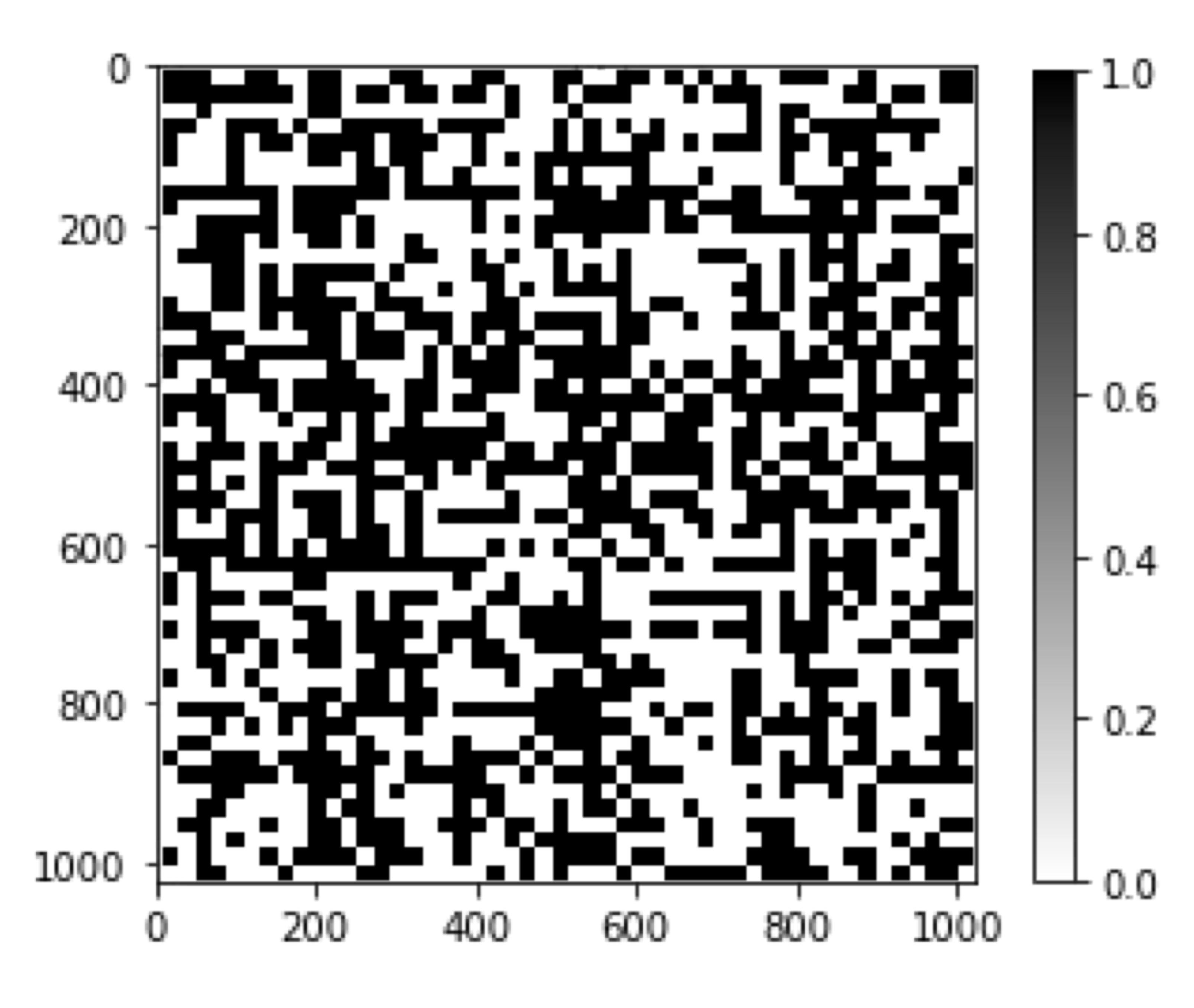}
\caption{Sinkhorn gets  $51\%$ sparsity.}
\label{sinkhorn-sparsity-color-transfer}
\end{figure}
\end{minipage}

\subsection{Approximation Error between UOT and OT}
\subsubsection{Approximation of OT Transportation Plan and OT Retrieval Gap }
\begin{minipage}[t]{0.3\textwidth}
\begin{figure}[H]
\includegraphics[width=1\textwidth]{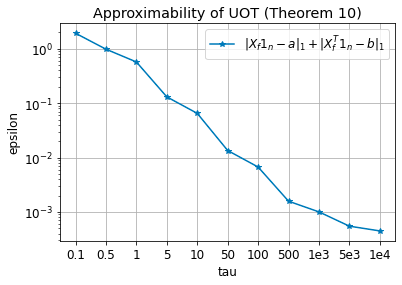}
\caption{Approximability of OT via UOT  (Theorem \ref{OT_solution_approx}) on synthetic data. The gap approaches $0$ with rate $O\left( \frac{n}{\tau} \right)$, as $\tau$ grows.} 
\label{fig_010}
\end{figure}
\end{minipage} 
\hfill
\begin{minipage}[t]{0.3\textwidth}
\begin{figure}[H]
\includegraphics[width=1\textwidth]{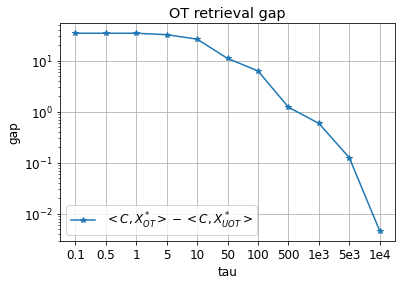}
\caption{OT retrieval of GEM-OT (Theorem \ref{gem_ot1}) on CIFAR-10. This empirically verifies vanishing OT retrieval gap for $\tau \to \infty$.}
\label{OT_retrieval_gap}
\end{figure}
\end{minipage}
\hfill
\begin{minipage}[t]{0.3\textwidth}
\begin{figure}[H]
\centering
\includegraphics[  width=0.95\textwidth]{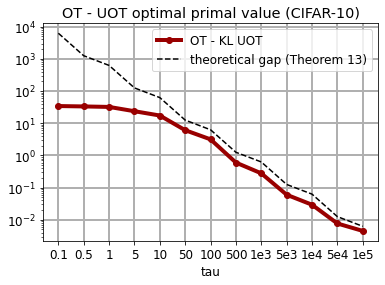}
\caption{The error $\OT(\aA, \bB)- \UOT_{\KL}(\aA, \bB)$ (Theorem \ref{thrm: distance_error}) in terms of the empirical gap (red) and the theoretical gap (black dotted) on the CIFAR-10. }
\label{fig_004}
\end{figure}
\end{minipage}
\vspace{.07in}

We investigate the approximability of OT via UOT in terms of transportation plan (Theorem \ref{OT_solution_approx}) on synthetic data where $\alpha = \beta = 1$ and entries of $C$ are drawn uniformly from $[0.5, 1]$. The result is depicted in Figure \ref{fig_010}.
Next, we verify the error incurred by GEM-OT (Algorithm \ref{alg: gem-ot}) for retrieving OT solution from UOT. The choice of $\tau$ in GEM-OT equivalently means that the retrieval gap $\varepsilon $ established in Theorem \ref{gem_ot1} can be upperbounded by $\frac{16 \|C\|_\infty n \gamma}{\tau}$ where $\gamma = \|C\|_\infty + \varepsilon/32$. To assert such vanishing effect of the error as $\tau$ grows large, we test on CIFAR-10 dataset and report the OT retrieval gap in Figure \ref{OT_retrieval_gap}. We note that  the OT retrieval gap  also upperbounds the approximation error in transportation plan, so Figure \ref{OT_retrieval_gap} also serves as empirical verification of Theorem \ref{OT_solution_approx} on CIFAR-10 data. \\

\subsubsection{Approximation Error between UOT and OT Distances}

On CIFAR-10 data, we measure the gap $\OT(\aA, \bB)-\UOT_{\KL}(\aA, \bB)$ of transport distance between OT and UOT (Theorem \ref{thrm: distance_error}) in Figure \ref{fig_004}. The result illustrates  the sharpness  of our theoretical bound and the trade-off between $\tau$ and approximation error, which can further elucidate more refined tuning of $\tau$. 

\section{Conclusion}\label{sec_conclusion}

We have developed a new algorithm for solving the discrete UOT problem using gradient-based methods. Our method achieves the complexity of $\widetilde{O}\big( \kappa n^2 \big) \big)$, where $\kappa$ is the condition number depending on only the two input measures. This is significantly better than the only known complexity ${O}\big(\frac{\tau n^2 \log(n)}{\varepsilon}  \log \big(\frac{\log(n)}{{\varepsilon}}\big)\big)$ of Sinkhorn algorithm \citep{pham2020unbalanced}  to the  UOT problem  in terms of $\varepsilon$ and $\tau$ dependency. In addition, we are the first to theoretically analyze the approximability of the UOT problem to find an OT solution. Our numerical results show the efficiency of the new method and the tightness of our theoretical bounds. We believe that our analysis framework can be extended to study the convergence results of stochastic methods, widely used in ML applications for handling large-scale datasets \citep{pmlr-v80-nguyen18c-hogwild} or under nonlinear stochastic control settings \citep{hoang_nonlinear_sa}. Our results on approximation error can open up directions that utilize UOT to approximate standard OT and even partial OT \citep{nguyen2023partial-ot-apdagd} using first-order methods.




\appendix


\section*{APPENDIX}

\section{Sinkhorn Algorithm}
\label{sinkhorn_sucks}

\subsection{Regularity Conditions} 
\label{sinkhorn_assumption}
In this Section, we will show that the assumptions used by the Sinkhorn \citep{pham2020unbalanced} algorithm are equivalent to our stated regularity conditions (S1)-(S4). 
(S4) is clearly stated in their list of regularity conditions. To claim their complexity, in  \citep[Corollary1]{pham2020unbalanced}, they further require the condition '$R = O(\frac{1}{\eta} \|C\|_\infty)$'. Following their discussion in  \citep[Section 3.1]{pham2020unbalanced}, a necessary condition for it to hold is $\norm{\log(\aA)}_{\infty} = O(\log(n))$ (and resp. $\norm{\log(\bB)}_{\infty} = O(\log(n))$). To see that $\norm{\log(\aA)}_{\infty} = O(\log(n))$ (and resp. $\norm{\log(\bB)}_{\infty} = O(\log(n))$) is equivalent to (S1)-(S3), we note that if $a_{min} = 0$ then $\norm{\log(\aA)}_{\infty} = \infty$ (contradiction!), and $\norm{\log(\aA)}_{\infty} = \max\{|\log(a_{min})|, |\log(a_{max})| \}$.


\subsection{The dependency of Sinkhorn complexity on 
  $\log(\alpha+\beta)$}


We note that in its pure form, Sinkhorn's complexity has $\log(\alpha+\beta)$, which is excluded in \citep{pham2020unbalanced} under their regularity condition (S4) that $\alpha, \beta$ are constants. To see this, note that in the proof of \citep[Corollary1]{pham2020unbalanced}, their quantity $U$ has the complexity $ O((\alpha+\beta)\log(n))$ and appears in Sinkhorn's final complexity as $\log(U)$.

\section{Supplementary Lemmas and Theorems}
\label{appen_supple}

\begin{lemma}
\label{max_smooth}
For $x, y \in \Br$, we have: $|\max\{0, x\} - \max\{0, y\}| \leq |x - y|$
\end{lemma}
\begin{proof}
We have: 
\begin{align*}
    \max\{0, x\} = \max\{0, x-y + y\} &\leq \max\{0, y\} + \max\{0, x-y\} \leq\max\{0, y\} + |x-y|  \\
    \max\{0, x\} - \max\{0, y\} &\leq |x-y|
\end{align*}
By symmetry, $ \max\{0, y\} - \max\{0, x\} \leq |x-y|$. Thus, $|\max\{0, x\} - \max\{0, y\}| \leq |x - y|$.
\end{proof}

\begin{theorem}[\cite{frenchel}]
\label{thm:frenchel_dual}
Let $(E, E^*)$ and $(F, F^*)$ be two couples of topologically paired spaces. Let $A : E \rightarrow F$ be a continuous linear operator and $A^* : F^* \rightarrow E^*$ its adjoint. Let f and g be lower semicontinuous and proper convex functions defined on E and F respectively. 
If there exists $x \in dom f$ such that g is continuous at $Ax$, then 
\begin{align*}
    \sup_{x \in E} -f(-x) - g(Ax)= \inf_{y^* \in F^*} f^*(A^* y^*) + g^*(y^*)
\end{align*}
Moreover, if there exists a maximizer $x \in E$ then there exists $y^* \in F^*$ satisfying $Ax \in \partial g^*(y^*)$ and $A^* y^* \in \partial f(-x)$.
\end{theorem}

\begin{lemma}[Lemma 7, \citep{altschuler2017near}]
\label{proj_algo}
For the inputs $X\in \Br^{n\times n}, \aA\in \Br^n$ and $\bB\in \Br^n$,  Algorithm \ref{alg: proj} takes $O(n^2)$ time to output a matrix $Y\in \Pi(\aA, \bB)$ satisfying:
\begin{align}
    \|Y-X\|_1 \leq 2 [\|X \one_n - \aA \|_1 + \|X^{\top} \one_n - \bB \|_1]
\end{align}
\end{lemma}



\section{Missing Proofs}

\subsection{Proof of Lemma \ref{lma: acce_proximal_guarantee1}}
\label{appen_acce_proximal_guarantee1}

Let $\xX^* =  (\uU^*, \vV^*, \tT^*)$ be an optimal solution to \eqref{dual_prob1}.
Because of the projection step on line \ref{projection_step1} of Algorithm \ref{alg: acc-proximal-gd-uot}, we have $\underline{\xX}^t, \xX^t \in V_D $ for $ t=0,\dots, k$. By Lemma \ref{sparse_UOT_lipschitz_lemma}, $f_\eta(\xX)$ is   $L$-smooth and $w_\eta(\xX)$ is $\mu$-strongly convex on $
conv\{\xX^*, \underline{\xX}^0, \xX^0, \underline{\xX}^1, \xX^1, ...\}$
The GEM-UOT is based on the GEM being applied to solve the dual objective \eqref{equi_dual_obj1}, which, from \citep[Theorem 2.1]{Lan2018RandomGE} and the proof therein, thus guarantees that the dual iterate $\underline{\xX}^k$ converges to the dual solution $\xX^*$ of \eqref{equi_dual_obj1} by:
\begin{align}
\| \underline{\xX}^k - \xX^*\|_2^2 &\leq \frac{4 \Delta_{0, \sigma_0} k \zeta^k (1-\zeta)}{\mu(1-\zeta^k)}. \label{baby1.1}
\end{align}
Furthermore, we have:
\begin{align}
\nonumber
     \frac{ k \zeta^k (1-\zeta)}{1-\zeta^k} &= \big( \sum_{t=1}^k \frac{\zeta^t}{\zeta^t}  \big)  \frac{ \zeta^k (1-\zeta)}{1-\zeta^k} \leq  \big( \sum_{t=1}^k \frac{\zeta^t}{\zeta^{1.5 t}}  \big)  \frac{ \zeta^k (1-\zeta)}{1-\zeta^k} \\
     \nonumber
     &\leq \frac{1-\zeta^{k/2}}{\zeta^{k/2}(1-\zeta^{1/2})}  \frac{ \zeta^k (1-\zeta)}{1-\zeta^k}\\
     &\leq 2\zeta^{k/2} \label{baby1.2}
\end{align}
Now, plugging \eqref{baby1.2} into \eqref{baby1.1} gives:
\begin{align}
    \| \underline{\xX}^k - \xX^*\|_2^2 &\leq \frac{8 \Delta_{0, \sigma_0} \zeta^{k/2}}{\mu}     \label{baby1}
\end{align}
Line \ref{convert_dual_primal} of Algorithm \ref{alg: acc-proximal-gd-uot} further converts the dual iterate $\underline{\xX}^k$ into the primal iterate $X^k$ as the final output. Next, we show that the convergence of the dual iterate in the sense of \eqref{baby1} leads to the convergence of the primal iterate $X^k$ to the primal solution $X^\eta$ of the primal objective \eqref{sparse_UOT_func}.
From the definition of $X^k$ in Algorithm \ref{alg: acc-proximal-gd-uot} and in view of Corollary \ref{duality_theorem},  we have:
\begin{align*}
    \|X^k - X^\eta\|_1 &= \frac{1}{2\eta} \sum_{i,j=1}^n |\underline{t}^k_{ij} - t^*_{ij} |\\
    &\leq  \frac{1}{2\eta} \| \underline{\xX}^k - \xX^* \|_1 \leq \frac{\sqrt{n^2+2n}}{2\eta} \| \underline{\xX}^k - \xX^* \|_2 \\
    &\overset{\eqref{baby1}}{\leq} \frac{\sqrt{2}n}{2\eta  \mu}  2\sqrt{2} \Delta_{0, \sigma_0}^{1/2}  \zeta^{k/4}  = \frac{2n\Delta_{0, \sigma_0}^{1/2}}{\eta  \mu}  (1- \frac{1}{1+\sqrt{1+16L/\mu}})^{k/4}\\
    &\leq \frac{2n\Delta_{0, \sigma_0}^{1/2}}{\eta  \mu}  \cdot exp\{- \frac{k}{4(1+\sqrt{1+16L/\mu})} \}.
\end{align*}
The condition that $\frac{2n\Delta_{0, \sigma_0}^{1/2}}{\eta  \mu}  \cdot exp\{- \frac{k}{4(1+\sqrt{1+16L/\mu})} \} \leq \min\{\frac{\varepsilon}{2 L_1} , \frac{\min\{a_{min}, b_{min}\} e^{-D/\tau}}{2}, \frac{\alpha+\beta}{2} \}$ is equivalent to: 
\begin{align*}
    k\geq 4( 1+\sqrt{1+16L/\mu} ) \log\bigg(\frac{4n\Delta_{0, \sigma_0}^{1/2} }{ \eta  \mu } \max\bigg\{\frac{L_1}{\varepsilon}, \frac{e^{D/\tau}}{\min\{a_{min}, b_{min}\}}, (\alpha+\beta)^{-1} \bigg\}\bigg) = K_0. 
\end{align*}
Therefore,
\begin{align}
\label{misc_bound1}
    \|X^k - X^\eta\|_1 \leq  \min\bigg\{\frac{\varepsilon}{2 L_1} , \frac{\min\{a_{min}, b_{min}\} e^{-D/\tau}}{2}, \frac{\alpha+\beta}{2} \bigg\}.
\end{align}


\subsection{Proof of Lemma \ref{baby_lemma1}}
\label{appen_baby_lemma1}
For $p= \frac{1}{2}\min\{a_{min}, b_{min}\} e^{-\frac{D}{\tau}}$ and $q = \alpha + \beta$ given in Section \ref{sec_assumptions}, we know that $ X^\eta\in U_{p,q}$ by Lemma \ref{opt_sol_bounds}. 
We now proceed to prove that $X^k \in U_{p,q}$. Consider any $i\in[1,n]$:  
\begin{align*}
    \sum_{j=1}^n X^k_{ij} &\geq \sum_{j=1}^n X^\eta_{ij} - |\sum_{j=1}^n (X^k_{ij} - X^\eta_{ij} )| \geq \sum_{j=1}^n X^\eta_{ij} - \|X^k - X^\eta \|_1 \overset{\eqref{gradient_cond1}}{=}  a_i e^{-u_i/\tau} - \|X^k - X^\eta \|_1 \\
    &\overset{\eqref{opt_sol_bounds3}}{\geq} a_{min} e^{-D/\tau} - \|X^k - X^\eta \|_1 \overset{\eqref{misc_bound1}}{\geq} \frac{a_{min} e^{-D/\tau}}{2} \geq p,\\
     \sum_{j=1}^n X^k_{ij} &\leq \sum_{j=1}^n X^\eta_{ij} + |\sum_{j=1}^n (X^k_{ij} - X^\eta_{ij} )| \leq \| X^\eta\|_1 + \|X^k - X^\eta \|_1\overset{\eqref{opt_sol_bounds1}, \eqref{misc_bound1}}{\leq} \frac{\alpha+\beta}{2}+ \frac{\alpha+\beta}{2} = \alpha + \beta. 
\end{align*}
Similarly, for any $j\in[1,n]$, we have: $p\leq \sum_{i=1}^n X^k_{ij}\leq q$, 
and thus conclude that $X^k\in U_{p,q}$.




\subsection{Proof of Lemma \ref{sparse_UOT_lipschitz_lemma}}
\label{lma10proof}

The gradient of $f_\eta(\xX = (\uU, \vV, \tT))$ can be computed as:
\begin{align*}
    &\frac{\partial f_\eta}{\partial u_i} = -a_i e^{-u_i/\tau}  -  \frac{\min\{a_{min}, b_{min}\}}{\tau}  e^{-D/\tau} u_i,\\
    &\frac{\partial f_\eta}{\partial v_j} = -b_j e^{-v_j/\tau}  -  \frac{\min\{a_{min}, b_{min}\}}{\tau}  e^{-D/\tau} v_j, \\
    &\frac{\partial f_\eta}{\partial t_{ij}} = 0. 
\end{align*}
The Hessian of $f_\eta(\xX = (\uU, \vV, \tT))$ can be computed as:
\begin{align*}
    &\frac{\partial^2 f_\eta}{\partial u_i \partial v_j} = \frac{\partial^2 f_\eta}{\partial u_i \partial t_{kj}} = \frac{\partial^2 f_\eta}{\partial v_i \partial t_{kj}} = 0, \\
    &\frac{\partial^2 f_\eta}{\partial u_i^2} = \frac{a_i}{\tau} e^{-u_i/\tau} -\frac{\min\{a_{min}, b_{min}\}}{\tau}  e^{-D/\tau},\\
    &\frac{\partial^2 f_\eta}{\partial v_j^2} = \frac{b_j}{\tau} e^{-v_j/\tau} -\frac{\min\{a_{min}, b_{min}\}}{\tau}  e^{-D/\tau}.
\end{align*}
If $\xX \in V_D$, then $u_i, v_j \leq D$ for all $i,j \in [n]$ which implies $\frac{\partial^2 f_\eta}{\partial u_i^2}, \frac{\partial^2 f_\eta}{\partial v_j^2}  \geq 0$ for all $i,j \in [n]$ and thereby  $\nabla^2 f_\eta(\xX) \succcurlyeq 0$ on $V_D$. We can then conclude that $f_\eta(\xX)$ is convex on $V_D$.

Now let us consider any $ \xX = (\uU, \vV, \tT),\xX' = (\uU', \vV', \tT') \in V_D$. By Mean Value Theorem, $\exists c_i \in range(u_i, u_i')$, such that $e^{-u_i/\tau} - e^{-u_i'/\tau} = -\frac{1}{\tau} e^{-c_i/\tau} (u_i - u_i')$ and $\exists d_j \in range(v_j, v_j')$, such that $e^{-v_j/\tau} - e^{-v_j'/\tau} = -\frac{1}{\tau} e^{-d_i/\tau} (v_j - v_j')$. We then obtain:
\allowdisplaybreaks
\begin{align*}
     \| \nabla f_\eta(\xX) - \nabla f_\eta(\xX') \|_2^2 =& \sum_{i=1}^n  \bigg[\big(\frac{a_i}{\tau} e^{-c_i/\tau}-\frac{\min\{a_{min}, b_{min}\}}{\tau}  e^{-D/\tau}\big) (u_i - u'_i)\bigg]^2 \\
     &+\sum_{j=1}^n  \bigg[\big(\frac{b_j}{\tau} e^{-d_j/\tau}-\frac{\min\{a_{min}, b_{min}\}}{\tau}  e^{-D/\tau}\big) (v_j - v'_j)\bigg]^2  \\
     \leq&  \sum_{i=1}^n  \big(\frac{a_i}{\tau} e^{-c_i/\tau}+\frac{\min\{a_{min}, b_{min}\}}{\tau}  e^{-D/\tau}\big)^2 (u_i - u'_i)^2 \\
     &+\sum_{j=1}^n  \big(\frac{b_j}{\tau} e^{-d_j/\tau}+\frac{\min\{a_{min}, b_{min}\}}{\tau}  e^{-D/\tau}\big)^2 (v_j - v'_j)^2  \\
     \leq&  \sum_{i=1}^n  \big(\frac{a_i}{\tau}  \frac{\alpha+\beta}{2a_i}+\frac{\min\{a_{min}, b_{min}\}}{\tau}  e^{-D/\tau}\big)^2 (u_i - u'_i)^2 \\
     &+\sum_{j=1}^n  \big(\frac{b_j}{\tau} \frac{\alpha+\beta}{2b_j}+\frac{\min\{a_{min}, b_{min}\}}{\tau}  e^{-D/\tau}\big)^2 (v_j - v'_j)^2  \\
     &\text{ (since $c_i\in range(u_i, u_i') $, $c_i \geq  \tau \log(\frac{2a_i}{\alpha+\beta}) $. Similarly, $d_j\geq  \tau \log(\frac{2b_j}{\alpha+\beta}) $)}\\
     \leq& \bigg(\frac{\alpha+\beta}{2\tau}+\frac{\min\{a_{min}, b_{min}\}}{\tau}  e^{-D/\tau}\bigg)^2 \|\xX - \xX'\|_2^2. \\
      \| \nabla f_\eta(\xX) - \nabla f_\eta(\xX') \|_2 \leq& \bigg(\frac{\alpha+\beta}{2\tau}+\frac{\min\{a_{min}, b_{min}\}}{\tau}  e^{-D/\tau}\bigg) \|\xX - \xX'\|_2,
\end{align*}
which implies that $f_\eta(\xX)$ is $L$-smooth with $L = \frac{\alpha+\beta}{2\tau}+\frac{\min\{a_{min}, b_{min}\}}{\tau}  e^{-D/\tau}$. 



Now let us consider $w_\eta(\xX)$. The gradient $w_\eta(\xX)$ can be computed as:
\begin{align*}
    \frac{\partial w}{\partial u_i} = \frac{\min\{a_{min}, b_{min}\}}{\tau}  e^{-D/\tau} u_i, \quad \frac{\partial w}{\partial v_j} = \frac{\min\{a_{min}, b_{min}\}}{\tau}  e^{-D/\tau} v_j, \quad \frac{\partial w}{\partial t_{ij}} &= \frac{1}{2\eta} t_{ij}.
\end{align*}
For any $ \xX = (\uU, \vV, \tT),\xX' = (\uU', \vV', \tT') \in V_D$, we have:
\begin{align*}
\left\langle \nabla w_\eta(\xX) - \nabla w_\eta(\xX') , \xX - \xX' \right\rangle &= \sum_{i=1}^n  \frac{\min\{a_{min}, b_{min}\}}{\tau}  e^{-D/\tau}  (u_i - u_i')^2  \\
&+ \sum_{j=1}^n  \frac{\min\{a_{min}, b_{min}\}}{\tau}  e^{-D/\tau}  (v_j - v_j')^2  + \sum_{i,j=1}^n \frac{1}{2\eta} (t_{ij} - t'_{ij})^2\\
&\geq \min\bigg\{ \frac{\min\{a_{min}, b_{min}\}}{\tau}  e^{-D/\tau} , \frac{1}{2\eta} \bigg\} \|\xX -\xX'\|_2^2.
\end{align*}
Therefore, $w_\eta(\xX)$ is $\mu$-strongly convex on $V_D$ with $\mu=  \min\bigg\{  \frac{\min\{a_{min}, b_{min}\}}{\tau}  e^{-D/\tau} , \frac{1}{2\eta} \bigg\}$.

\subsection{Projection Step of GEM-UOT}
\label{appen_projection_gemuot}
We describe the following projection step (i.e. line \ref{projection_step1} in  Algorithm \ref{alg: acc-proximal-gd-uot}) in more details:
\begin{align}
 \xX^t &= \mathcal{M}_{V_D\cap \Xx}(\widetilde{\yY}^t, \xX^{t-1}, \rho) = \argminE_{\xX = (\uU,\vV, \tT) \in V_D\cap \Xx} \{ G(\xX) := \langle \widetilde{\yY}^t,\xX \rangle + w_\eta(\xX) + \rho P( \xX^{t-1},\xX) \} \label{projection_inter1}
\end{align}
where we further expand the objective $G(\xX)$ as:
\begin{align}
\nonumber
    G(\xX) &= \langle \widetilde{\yY}^t,\xX \rangle + w_\eta(\xX) + \frac{\rho}{\mu} [ w_\eta(\xX) - w_\eta( \xX^{t-1}) - \langle \nabla w_\eta( \xX^{t-1}), \xX- \xX^{t-1} \rangle ]\\
    &= \big(1 + \frac{\rho}{\mu}  \big) w_\eta(\xX) + \langle \widetilde{\yY}^t - \frac{\rho}{\mu} \nabla w_\eta( \xX^{t-1}),\xX \rangle + \frac{\rho}{\mu}  \langle \nabla w_\eta( \xX^{t-1}), \xX^{t-1} \rangle.  \label{projection_inter2}
\end{align}
For abbreviation, we let $\widetilde{\yY}^t - \frac{\rho}{\mu} \nabla w_\eta( \xX^{t-1}) =  (\uU^c, \vV^c, \tT^c) \in \Br^{n^2+2n}$ where $\uU^c\in \Br^{n}, \vV^c\in \Br^{n}$ and $\tT^c \in \Br^{n^2}$ respectively are coefficients of the variables $\uU, \vV, \tT$ in the linear term of \eqref{projection_inter2}. We also define the constants $q_1 = \big(1+\frac{\rho}{\mu}  \big) \cdot \frac{\min\{a_{min}, b_{min}\}}{2\tau}  e^{-D/\tau} $ and $q_2 = \big(1+\frac{\rho}{\mu}  \big)/(4\eta)$.
The objective $G(\xX)$ can be written as the sum of quadratic functions as follows:
\begin{align}
    G(\xX) = ( q_1 \|\uU\|_2^2 + \langle \uU^c, \uU\rangle )+ ( q_1 \|\vV\|_2^2 + \langle \vV^c, \vV\rangle) + (q_2 \|\tT\|_2^2 +  \langle \tT^c, \tT \rangle ) + \frac{\rho}{\mu}  \langle \nabla w_\eta( \xX^{t-1}), \xX^{t-1} \rangle.  
\end{align}
And the projection step is thus equivalent to solving:
\begin{align}
    \min_{\xX = (\uU,\vV, \tT) \in V_D\cap \Xx} &G(\xX)  \label{projection_obj1} \\
    \nonumber
   \text{where } V_D\cap\Xx = \{(\uU, \vV, \tT) | &\uU, \vV \in \Br^{n}, \tT \in \Br^{n\times n}: \forall i, j \in[1,n], \\
   \nonumber
   & \tau \log(\frac{2a_i}{\alpha+\beta}) \leq u_i \leq D \\
   \nonumber
   &\tau \log(\frac{2b_j}{\alpha+\beta}) \leq  v_j \leq D \\
   \nonumber
   &0 \leq t_{ij},\\
   \nonumber
   &u_i + v_j - C_{ij} \leq t_{ij}\}.
\end{align}
Next, we aim to derive a simplified  equivalence form of \eqref{projection_obj1} that is obtained by expressing the optimal $\tT$ with respect to fixed $\uU, \vV$, and consequently optimizing over $\uU, \vV$. In particular, we consider the mapping $\tT(\uU, \vV): \Br^{2n} \to \Br^{n^2}$ defined by:
\begin{align}
    t_{ij}(\uU, \vV) =  \max\big\{0, -\frac{t_{ij}^c}{q_2}, u_i + v_j - C_{ij}   \big\}, \quad \forall i, j \in [n], \label{conditioned_t1}
\end{align}
the restricted domain $\widetilde{V}_D$ as follows:
\begin{align}
\nonumber
    \widetilde{V}_D = \{(\uU, \vV) | &\uU, \vV \in \Br^{n}: \forall i \in[n], \\
   \nonumber
   & \tau \log(\frac{2a_i}{\alpha+\beta}) \leq u_i \leq D \\
   &\tau \log(\frac{2b_i}{\alpha+\beta}) \leq  v_i \leq D\}, \label{newVD}
\end{align}
and the objective:
\begin{align}
    \widetilde{G}(\uU, \vV) = G(\uU, \vV, \tT(\uU, \vV)). \label{newprojection}
\end{align}
The next Lemma establishes the equivalence form for the projection problem \eqref{projection_obj1}.
\begin{lemma}
 \label{lma:equivalence_project}
 The following holds:
 \begin{align}
      \min_{\xX = (\uU,\vV, \tT) \in V_D\cap \Xx} &G(\xX)  =  \min_{(\uU,\vV) \in \widetilde{V}_D} \widetilde{G}(\uU, \vV) 
 \end{align}
\end{lemma}
\begin{proof}
For any $(\uU, \vV) \in \widetilde{V}_D$, we have $\tT(\uU, \vV)$ that is defined as in \eqref{conditioned_t1} satisfies:
\begin{align*}
    0\leq t_{ij}(\uU, \vV),\\
    u_i+v_j-C_{ij} \leq t_{ij}(\uU, \vV),
\end{align*}
Thus, we have $(\uU, \vV, \tT(\uU, \vV)) \in V_D$, which implies:
\begin{align}
\nonumber
   &\min_{\xX = (\uU,\vV, \tT) \in V_D\cap \Xx} G(\xX) \leq  G(\uU, \vV, \tT(\uU, \vV)) = \widetilde{G}(\uU, \vV), \quad \forall (\uU, \vV) \in \widetilde{V}_D
    \\
   \therefore &\min_{\xX = (\uU,\vV, \tT) \in V_D\cap \Xx} G(\xX) \leq \min_{(\uU,\vV) \in \widetilde{V}_D} \widetilde{G}(\uU, \vV).  \label{equivalence_project_baby1}
\end{align}
Now, take any $\xX = (\uU, \vV, \tT) \in V_D \cap \Xx$. We fix $\uU$ and $\vV$, and optimize over all $\tT'$ such that $(\uU, \vV, \tT') \in V_D \cap \Xx$, i.e. those $\tT'$ with $0 \leq t'_{ij}$ and $u_i + v_j - C_{ij} \leq t'_{ij}$. Then, we have:
\begin{align}
    G(\xX  = (\uU, \vV, \tT)) \geq  \min_{\tT': (\uU, \vV, \tT') \in V_D \cap \Xx} G(\uU, \vV, \tT').  \label{equivalence_project_baby2}
\end{align}
Note that the optimization problem over $\tT'$ in the RHS of \eqref{equivalence_project_baby2} has separable quadratic objective and separable linear constraints.  Thus, simple algebra gives $\tT(\uU, \vV)$ being defined in \eqref{conditioned_t1} as the optimal solution to such problem, i.e.
\begin{align}
    \tT(\uU, \vV) = \argminE_{\tT': (\uU, \vV, \tT') \in V_D \cap \Xx} G(\uU, \vV, \tT').  \label{equivalence_project_baby3}
\end{align}
From \eqref{equivalence_project_baby2} and \eqref{equivalence_project_baby3}, we have:
\begin{align}
\nonumber
      G(\xX  = (\uU, \vV, \tT)) &\geq  G(\uU, \vV, \tT(\uU, \vV)) = \widetilde{G}(\uU, \vV)\\
      \nonumber
      &\geq \min_{(\uU,\vV) \in \widetilde{V}_D} \widetilde{G}(\uU, \vV)\\
      \therefore \min_{\xX = (\uU,\vV, \tT) \in V_D\cap \Xx} G(\xX)   &\geq \min_{(\uU,\vV) \in \widetilde{V}_D} \widetilde{G}(\uU, \vV),\label{equivalence_project_baby4}
\end{align}
where the second inequality is because any $\xX = (\uU, \vV, \tT) \in V_D \cap \Xx$ implies $(\uU, \vV) \in \widetilde{V}_D$ by definition. Finally, combining \eqref{equivalence_project_baby1} and \eqref{equivalence_project_baby4} gives the desired statement.
\end{proof}

Since the  equivalence form  $\min_{(\uU,\vV) \in \widetilde{V}_D} \widetilde{G}(\uU, \vV) $  of the projection step is strongly convex smooth optimization over box constraints, we can adopt projected quasi-Newton method  with superlinear convergence rate \citep[Theorem 11.2]{projectedNewton}, thereby resulting in $\widetilde{O}(1)$ iteration complexity. 
 Furthermore, the per iteration cost is dominated by the gradient and Hessian computation of $\widetilde{G}$, which can be achieved efficiently in $O(n^2)$ as follows. In particular,  we have for $\forall i, j \in [n]$:
\begin{align*}
    \frac{\partial \widetilde{G}}{\partial u_i} &= 2q_1 u_i + u^c_i + 2 \sum_{j=1}^n \max\bigg\{0 , \frac{t_{ij}^c}{q_2}, u_i +v_j - C_{ij} +\frac{t_{ij}^c}{q_2}  \bigg\},\\
    \frac{\partial \widetilde{G}}{\partial v_j} &= 2q_1 v_j + v^c_j + 2 \sum_{i=1}^n \max\bigg\{0 , \frac{t_{ij}^c}{q_2}, u_i +v_j - C_{ij} +\frac{t_{ij}^c}{q_2}  \bigg\},\\
    \frac{\partial^2 \widetilde{G}}{\partial u_i \partial u_j} &= \begin{cases} 2q_1 + 2 \sum_{k=1}^n \mathbbm{1}\big\{u_i +v_k - C_{ik}  \geq \max\{0, -\frac{t_{ik}^c}{q_2}\}   \big\},  \quad \text{ if $i= j$} \\
    0, \quad \text{ otherwise}
    \end{cases},\\
    \frac{\partial^2 \widetilde{G}}{\partial v_i \partial v_j} &= \begin{cases} 2q_1 + 2 \sum_{k=1}^n \mathbbm{1}\big\{u_k +v_j - C_{kj}  \geq \max\{0, -\frac{t_{kj}^c}{q_2}\}   \big\},  \quad \text{ if $i= j$} \\
    0, \quad \text{ otherwise}
    \end{cases},\\
    \frac{\partial^2 \widetilde{G}}{\partial u_i \partial v_j} &= 2 \mathbbm{1}\big\{u_i +v_j - C_{ij}  \geq \max\{0, -\frac{t_{ij}^c}{q_2}\}   \big\}.
\end{align*}
Thus, we conclude that the projection step can be solved in  $\widetilde{O}(n^2)$ time.







\subsection{Proof of Lemma \ref{rgem_convex_lipschitz_lemma}}
\label{appen_lma16}

The gradient of $h_a(\xX = (\uU, \vV))$ can be computed as:
\begin{align*}
    &\frac{\partial h_a}{\partial u_i} = -a_i e^{-u_i/\tau} + \sum_{j = 1}^n\frac{\max\{0, u_i+v_j-C_{ij}\}}{2\eta},\\
    &\frac{\partial h_a}{\partial v_j} = -b_j e^{-v_j/\tau} + \sum_{i = 1}^n\frac{\max\{0, u_i+v_j-C_{ij}\}}{2\eta}. 
\end{align*}

Let us consider any $ \xX = (\uU, \vV),\xX' = (\uU', \vV') \in V_a$. By Mean Value Theorem, $\exists c_i \in range(u_i, u_i')$, such that $e^{-u_i/\tau} - e^{-u_i'/\tau} = -\frac{1}{\tau} e^{-c_i/\tau} (u_i - u_i')$ and $\exists d_j \in range(v_j, v_j')$, such that $e^{-v_j/\tau} - e^{-v_j'/\tau} = -\frac{1}{\tau} e^{-d_i/\tau} (v_j - v_j')$. We then obtain:
\begin{align*}
     &\| \nabla h_a(\xX) - \nabla h_a(\xX') \|_2^2 \leq 2\sum_{i=1}^n  \bigg[\frac{a_i}{\tau} e^{-c_i/\tau} (u_i - u'_i)\bigg]^2 + 2 \sum_{j=1}^n  \bigg[\frac{b_j}{\tau} e^{-d_j/\tau} (v_j - v'_j)\bigg]^2 \\
     &+\sum_{i=1}^n \frac{2}{4\eta^2} [\sum_{j=1}^n (u_i - u_i' + v_j-v_j')]^2 + \sum_{j=1}^n \frac{2}{4\eta^2} [\sum_{i=1}^n (u_i - u_i' + v_j-v_j')]^2\\
     &\text{(By Lemma \ref{max_smooth} and simple algebra $(x+y)^2\leq 2(x^2+y^2)$)} \\
     \leq&  \sum_{i=1}^n  2\left(\frac{a_i}{\tau} e^{-c_i/\tau}\right)^2 (u_i - u'_i)^2 +\sum_{j=1}^n  2\left(\frac{b_j}{\tau} e^{-d_j/\tau} \right)^2 (v_j - v'_j)^2 \\
     &+ \sum_{i=1}^n \frac{1}{\eta^2} \{\sum_{j=1}^n [(u_i - u_i')^2 + (v_j-v_j')^2]\} + \sum_{j=1}^n \frac{1}{\eta^2} \{\sum_{i=1}^n [(u_i - u_i')^2 + (v_j-v_j')^2]\} \\
     &\leq  \sum_{i=1}^n  2\left(\frac{a_i}{\tau}  \frac{\alpha+\beta}{2a_i}\right)^2 (u_i - u'_i)^2 +\sum_{j=1}^n  2 \big(\frac{b_j}{\tau} \frac{\alpha+\beta}{2b_j}\big)^2 (v_j - v'_j)^2 +\frac{2n}{\eta^2}\|\xX - \xX'\|_2^2\\
     &\text{ (since $c_i\in range(u_i, u_i') $, $c_i \geq  \tau \log(\frac{2a_i}{\alpha+\beta}) $. Similarly, $d_j\geq  \tau \log(\frac{2b_j}{\alpha+\beta}) $)}\\
     &= \left[\frac{1}{2}\cdot \bigg( \frac{\alpha+\beta}{\tau}\bigg)^2 + \frac{2n}{\eta^2} \right]\|\xX - \xX'\|_2^2. \\
      \therefore &\| \nabla h_a(\xX) - \nabla h_a(\xX') \|_2 \leq \bigg(\frac{\alpha+\beta}{\tau} + \frac{2\sqrt{n}}{\eta} \bigg)\cdot \|\xX - \xX'\|_2.
\end{align*}
This implies that $h_a(\xX)$ is $L$-smooth with $L_a = \frac{\alpha+\beta}{\tau} +  \frac{2\sqrt{n}}{\eta}$.

\section{Experimental Result for GEM-RUOT on Fashion-MNIST dataset}
\label{sec:mnist1}
\begin{figure}[H]
\centering
\includegraphics[width=0.6\textwidth]{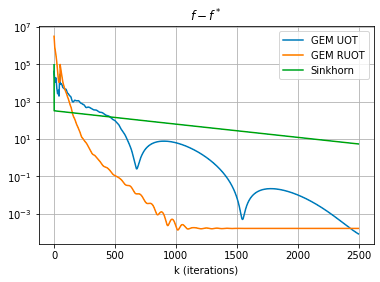}
\caption{Comparison of $f$ primal gap between the Sinkhorn algorithm and GEM algorithms for $1500$ iterations.}
\label{fig_009}
\end{figure}

We test GEM-UOT (Algorithm \ref{alg: acc-proximal-gd-uot}), GEM-RUOT (Algorithm \ref{alg: acc-proximal-gd-uot-convex}) and Sinkhorn on the Fashion-MNIST dataset using the same setting as the experiment on CIFAR dataset in Section \ref{sec:real_data_exp}  and report the result in Figure \ref{fig_009}. 



\vskip 0.2in
\bibliography{sample}

\end{document}